\documentclass[a4paper]{amsart}
\usepackage{latexsym}
\usepackage{mathtools}
\usepackage{amsmath}
\usepackage{amsfonts}
\usepackage{mathrsfs}
\usepackage{amscd}
\usepackage{amssymb}
\usepackage{amsfonts}
\usepackage{dsfont}
\usepackage{float}
\usepackage{esint}
\usepackage{caption}
\usepackage{hyperref}
\usepackage{tikz-cd}
\usepackage{slashed}
\usepackage[utf8]{inputenc}
\usepackage[T1]{fontenc}   
\usepackage{graphicx}
\usepackage{color}
\allowdisplaybreaks 
\definecolor{Blue}{rgb}{0.,0.,1.}
\definecolor{Red}{rgb}{1.,0.,0.}
\definecolor{Green}{rgb}{0.,1.,0.}

\let\origmaketitle\maketitle
\def\maketitle{
  \begingroup
	\origmaketitle
  \endgroup
	}	
	
\makeatletter
\renewcommand{\author}[2][]{%
  \def\@tempa{#1}
  \ifx\@empty\authors
    \ifx\@tempa\@empty
      \gdef\shortauthors{#2}%
    \else
      \gdef\shortauthors{#1}%
    \fi
    \gdef\authors{\author{#2}}%
  \else
    \ifx\@tempa\@empty
      \g@addto@macro\shortauthors{\and#2}%
    \else
      \g@addto@macro\shortauthors{\and#1}%
    \fi
    \g@addto@macro\authors{\and\author{#2}}%
  \fi
}
\renewcommand{\address}[2][]{\g@addto@macro\authors{\address{#1}{#2}}}
\def\@setauthors{%
  \begin{center}%
    \footnotesize
    \vspace{20pt}
    \let\and\@empty
    \def\author##1{\advance\@tempcnta\@ne}%
    \def\address##1##2{\advance\@tempcntb\@ne}%
    \@tempcnta=\z@  \@tempcntb=\z@
    \authors
    \ifnum\@tempcnta>\@ne \ifnum\@tempcntb=\@ne
        \oneaddress
      \else
        \sepaddresses
      \fi
    \else
      \oneaddress
    \fi
  \end{center}%
}
\def\oneaddress{%
  \begingroup
  \let\author\@iden \let\address\@gobbletwo
  \renewcommand{\andify}{%
    \nxandlist{\unskip, }{\unskip{} and~}{\unskip, and~}}%
  \uppercasenonmath\authors
  \andify\authors
  \authors
  \endgroup
  \begingroup \let\and\relax \let\author\@gobble
  \def\address##1##2{\unskip\\[10pt] \itshape##2}%
  \authors
  \endgroup
}
\def\sepaddresses{%
  \begingroup
    \baselineskip10\p@\relax
    \def\address##1##2{ ({\itshape##2}\/)}
    \def\author##1{\def\temp{##1}\leavevmode\uppercasenonmath\temp\temp}%
    \nxandlist
      {,\\[\baselineskip]}
      {\\[\baselineskip] \textsc{\lowercase{and}}\\[\baselineskip]}
      {,\\[\baselineskip]\textsc{\lowercase{and}}\\[\baselineskip]}
      \authors 
    \authors
  \endgroup
}
\def\maketitle{\par
  \@topnum\z@
  \@setcopyright
  \thispagestyle{firstpage}%
  \uppercasenonmath\shorttitle
  \ifx\@empty\shortauthors \let\shortauthors\shorttitle
  \else
    \newcommand{\@xuppercasenonmath}[1]{\toks@\@emptytoks
      \@xp\@skipmath\@xp\@empty##1$$%
      \edef##1{\@nx\protect\@nx\@upprep\the\toks@}}%
    \@xuppercasenonmath\shortauthors
    \def\@@and{AND}
    \renewcommand{\andify}{%
      \nxandlist{\unskip, }{\unskip{ }\@@and{ }}{\unskip, \@@and{ }}}%
    \andify\shortauthors
  \fi
  \@maketitle@hook
  \begingroup
  \@maketitle
  \endgroup
  \c@footnote\z@
  \@cleartopmattertags
}
\def\@maketitle{%
  \normalfont\normalsize
  \let\@makefntext\noindent
  \@adminfootnotes
  \ifx\@empty\addresses\else \@footnotetext{\@setotheraddresses}\fi
  \global\topskip68\p@\relax
  \@settitle
  \ifx\@empty\authors \else \@setauthors \fi
  \ifx\@empty\@dedicatory
  \else
    \baselineskip26\p@
    \vtop{\centering{\footnotesize\itshape\@dedicatory\@@par}%
      \global\dimen@i\prevdepth}\prevdepth\dimen@i
  \fi
  \toks@\@xp{\shortauthors}\@temptokena\@xp{\shorttitle}%
  \edef\@tempa{\@nx\markboth{\the\toks@}{\the\@temptokena}}\@tempa
  \@setabstract
  \normalsize
  \if@titlepage
    \newpage
  \else
    \dimen@34\p@ \advance\dimen@-\baselineskip
    \vskip\dimen@\relax
  \fi
} 
\renewcommand{\thanks}[1]{%
  \ifx\@empty\thankses
    \gdef\thankses{\thanks{#1}}%
  \else
    \g@addto@macro\thankses{\endgraf\thanks{#1}}%
  \fi}
\def\@setthanks{\def\thanks##1{\noindent##1\@addpunct.}\thankses}
\renewcommand{\curraddr}[2][]{%
  \ifx\@empty\addresses
    \gdef\addresses{\curraddr{#1}{#2}}%
  \else
    \g@addto@macro\addresses{\endgraf\curraddr{#1}{#2}}%
  \fi}
\renewcommand{\email}[2][]{%
  \ifx\@empty\addresses
    \gdef\addresses{\email{#1}{#2}}%
  \else
    \g@addto@macro\addresses{\endgraf\email{#1}{#2}}%
  \fi}
\renewcommand{\urladdr}[2][]{%
  \ifx\@empty\addresses
    \gdef\addresses{\urladdr{#1}{#2}}%
  \else
    \g@addto@macro\addresses{\endgraf\urladdr{#1}{#2}}%
  \fi}
\def\@setotheraddresses{%
  \def\curraddr##1##2{\noindent
    \emph{Current address\@ifnotempty{##1}{ of ##1}}:\space
      ##2\@addpunct.}%
  \def\email##1##2{\noindent
    \emph{E-mail address\@ifnotempty{##1}{ of ##1}}:\space
      \texttt{##2}}%
  \def\urladdr##1##2{\noindent
    \emph{WWW address\@ifnotempty{##1}{ of ##1}}:\space
      \texttt{##2}}%
  \addresses
}
\let\enddoc@text\relax
\makeatother

\newcounter{smallarabics}
\newenvironment{arabicenumerate}
{\begin{list}{{\normalfont\textrm{(\arabic{smallarabics})}}}
  {\usecounter{smallarabics}\setlength{\itemindent}{0cm}
   \setlength{\leftmargin}{5ex}\setlength{\labelwidth}{4ex}
   \setlength{\topsep}{0.75\parsep}\setlength{\partopsep}{0ex}
   \setlength{\itemsep}{0ex}}}
{\end{list}}

\newcounter{smallroman}

\newcommand{\ben}{\begin{arabicenumerate}}  
\newcommand{\een}{\end{arabicenumerate}}

\def\init{\setcounter{equation}{0}}

\newtheorem{theoreme}{Theorem }[section]
\newtheorem{proposition}[theoreme]{Proposition}

\newtheorem{lemma}[theoreme]{Lemma}
\newtheorem{definition}[theoreme]{Definition}

\newtheorem{remark}[theoreme]{Remark}
\newtheorem{example}[theoreme]{Example}
\newcommand{\beq}{\begin{equation}}
\newcommand{\eeq}{\end{equation}}

\newcommand{\bex}{\begin{example}}
\newcommand{\eex}{\end{example}}
\def\bel{\begin{lemma}}
\def\eel{\end{lemma}}
\def\bet{\begin{theoreme}}
\def\eet{\end{theoreme}}
\def\bed{\begin{definition}}
\def\eed{\end{definition}}
\def\ber{\begin{remark}}
\def\eer{\end{remark}}

\def\id{{\mathds{1}}}

\def\rr{{\mathbb R}}

\def\cc{{\mathbb C}}
\def\nn{{\mathbb N}}

\def\Re{{\rm Re}}

\def\bar{\overline}

\def\cinf{C^\infty}
\def\proof{
\noindent{\bf Proof.}\ \ }

\def\cY{{\mathcal Y}}

\def\cL{{\mathcal L}}
\def\cS{{\mathcal S}}

\def\cD{{\mathcal D}}
\def\cU{{\mathcal U}}

\def\cC{{\mathcal C}}
\def\cW{{\mathcal W}}

\def\i{{\rm i}}
\def\qed{$\Box$\medskip}
\def \p{ \partial}
\def\12{\frac{1}{2}}
\def\14{\frac{1}{4}}

\def\bbbone{{\mathchoice {\rm 1\mskip-4mu l} {\rm 1\mskip-4mu l}
{\rm 1\mskip-4.5mu l} {\rm 1\mskip-5mu l}}}
\def\one{\bbbone}
\def\cH{{\mathcal H}}

\def\coinf{C_0^\infty}

\def\cT{{\mathcal T}}
\def\cF{{\mathcal F}}
\def\cG{{\mathcal G}}
\def\cX{{\mathcal X}}

\def \p{ \partial}
\def\12{\frac{1}{2}}
\def\e{{\rm e}}

\def\Op{{\rm Op}}

\newcommand{\mat}[4]{\left(\begin{array}{cc}#1 &#2  \\ #3 &#4 \end{array}\right)}

\def\cE{{\mathcal E}}

\def\WF{{\rm WF}}
\makeatletter
\newcommand*{\defeq}{\mathrel{\rlap{%
                     \raisebox{0.3ex}{$\m@th\cdot$}}%
                     \raisebox{-0.3ex}{$\m@th\cdot$}}%
                     =}
\makeatother
\makeatletter
\newcommand*{\eqdef}{=\mathrel{\rlap{%
                     \raisebox{0.3ex}{$\m@th\cdot$}}%
                     \raisebox{-0.3ex}{$\m@th\cdot$}}%
                     }
\makeatother

\def\Sol{{\rm Sol}_{\rm sc}}

\def\cinfb{C^{\infty}_{\rm b}}
\def\rx{{\rm x}}

\DeclareMathOperator{\Ker}{Ker}
\DeclareMathOperator{\Dom}{Dom}

\def\dual{\!\cdot \!}
\def\CCR{{\rm CCR}}

\def\cB{\mathcal{B}}
\def\fA{{\mathfrak A}}

\def\zero{{\mskip-4mu{\rm\textit{o}}}}
\def\cC{{\mathcal C}}
\def\cN{{\mathcal N}}

\def\BT{{\rm  BT}}

\def\tosim{\xrightarrow{\sim}}

\def\zero{{\mskip-4mu{\rm\textit{o}}}}
\def\Cl{{\rm Cl}}
\def\SO{{\rm SO}}
\def\Spin{{\rm Spin }}
 
\def\CAR{{\rm CAR}}
\def\CARY{{\rm CAR(\cY,\nu)}}
	\def\maketitle{
  \begingroup
  \def\uppercasenonmath##1{} 
  \let\MakeUppercase\relax 
	\origmaketitle
  \endgroup
	}

\begin{document}
\pagestyle{plain}

\title{\large Hadamard states for quantized Dirac fields on Lorentzian manifolds of bounded geometry
}
\author{\normalsize Christian \textsc{G\'erard}}
\address{Universit\'e Paris-Saclay, D\'epartement de Math\'ematiques, 91405 Orsay Cedex, France}
\email{christian.gerard@math.u-psud.fr}
\date{June 2021}
\author{\normalsize Théo \textsc{Stoskopf} }
\email{theo.stoskopf@universite-paris-saclay.fr}
\keywords{Hadamard states, microlocal spectrum condition,  pseudo-differential calculus, Dirac equation, curved spacetimes}
\subjclass[2010]{ 81T20, 35S05, 35Q41,}
\begin{abstract} 
We consider Dirac equations on even dimensional Lorentzian manifolds of bounded geometry with a spin structure.  For the associated free quantum field theory, we construct pure Hadamard states using global pseudodifferential calculus on a Cauchy surface. We also give two constructions of Hadamard states for Dirac fields for arbitrary spacetimes with a spin structure.
\end{abstract}

\maketitle

\section{Introduction}\label{sec0}
\subsection{Quantum fields on curved spacetimes}\label{sec0.1}
In many situations of physical interest, for example the study of the early stages of the universe or of stellar collapse,  one is naturally led to the problem of constructing quantum field theories on a  curved spacetime. 

The symmetries of the Minkowski spacetime, which play such a fundamental role, are absent in curved spacetimes, except in some simple situations, like {\em stationary} or {\em static} spacetimes.
 Therefore, the traditional approach to quantum field theory has to be modified: one has first to perform 
 an {\em algebraic quantization}, which for free theories amounts to introducing an appropriate {\em phase space}, which is either a {\em symplectic} or an {\em Euclidean} space, in the {\em bosonic} or {\em fermionic} case.
 From such a phase space one can construct $\CCR$ or ${\rm CAR}$ $*${\em -algebras}, and actually 
 {\em nets} of $*$-algebras, each associated to a region of spacetime.
 
 The second step consists in singling out, among the many states on these $*$-algebras, the physically meaningful ones, which should resemble the Minkowski vacuum, at least in the vicinity of any point of the spacetime. 
 This leads to the notion of {\em Hadamard states}, which were originally defined by requiring that their two-point functions have a specific asymptotic expansion near the diagonal, called the {\em Hadamard expansion}.
 
 A very important progress was made by Radzikowski, \cite{R1, R2}, who introduced the characterization of Hadamard states by the {\em wavefront set} of their two-point functions. 
Numerous papers have been devoted to the study of linear scalar fields,  on curved spacetimes, but much less has been done for fields with higher spin, in particular for {\em Dirac fields}, to which the present paper is devoted.

\subsection{Hadamard states for Dirac fields}\label{sec0.2}
Let $(M, g)$ a spacetime, i.e. an orientable and time orientable Lorentzian manifold. We will always assume that $ \dim M$ is even. If $(M, g)$ admits a spin structure, one can canonically define  {\em Dirac operators}
\[
D= g^{\mu\nu} \gamma(e_{\mu})\nabla^{S}_{e_{\nu}}+ m
\]
acting on sections of the spinor bundle $S(M)$. Here $(e_{\mu})_{0\leq \mu\leq d}$ is a local frame of $TM$, $g^{\mu\nu}$ is the inverse metric, $\gamma(e_{\mu})$ the Clifford multiplication and $\nabla^{S}$ the spin connection. The mass of the Dirac field is described by a real function $m: M\to \rr$ or more generally by a map $m: M\to L(S(M))$ which is selfadjoint for the spinor scalar product.

If  $(M,g)$ is globally hyperbolic, one can define the {\em retarded/advanced inverses} $G_{\rm ret/adv}$  of $D$  and the {\em causal propagator} $G= G_{\rm ret}- G_{\rm adv}$. Using the causal propagator one can, by the classic paper of Dimock \cite{Di}, construct the $\CAR$ $*$-{\em algebra} $\CAR(D)$, which describes  a free quantum Dirac field on $M$.

A {\em state} $\omega$ on $\CAR(D)$ (more precisely a gauge invariant quasi-free state), is  completely specified  by fixing a pair of linear continuous maps 
\[
\Lambda^{\pm}: \coinf(M; S(M))\to\cD'(M; S(M))
\]
 called the {\em covariances} of $\omega$ which satisfy:
\[
\begin{array}{l}
	\Lambda^{\pm} \geq 0 \hbox{ for }(\cdot|\cdot)_{M}, \\[2mm]
	\Lambda^+ + \Lambda^- = \i G,\\[2mm]
	D \circ \Lambda^{\pm} = \Lambda^{\pm} \circ D = 0,
\end{array}	
\]
where the (non positive) scalar product $(\cdot| \cdot)_{M}$ is defined in \eqref{burk}.

{\em Hadamard states} for Dirac fields were originally defined in \cite{K,V}, by specifying the short distance behavior of $\Lambda^{\pm}$ near the diagonal, in an analogous way to the case of Klein-Gordon fields.  The microlocal definition of Hadamard states for Dirac fields  was first introduced by Hollands in \cite{Ho}. 

The microlocal definition of Hadamard states is analogous to the case of scalar fields, although a more refined description of the microlocal singularities of $\Lambda^{\pm}$, based on the notion of the {\em polarization set} due to Dencker \cite{De}, was given in \cite{Ho}. 
The equivalence of the two definitions, which for scalar fields was one of the important results of Radzikowski, was proved in \cite{Ho}, see also the work of Kratzert \cite{Kr}.

These results were extended to more general hyperbolic equations acting on sections of vector bundles by Sahlman and Verch \cite{SV}. Later on Sanders described Dirac fields on four dimensional spacetimes as a locally covariant quantum field theory using the language of category theory. He also filled several  gaps in the existing literature.

Another approach to Hadamard states  has been advocated by Finster, called the {\em fermionic projector method}, see e.g. \cite{FR1, FR2}. 
Its scope seems for the moment  limited to rather special classes of spacetimes. Nevertheless its relationship with the pseudodifferential approach that we will use in this paper is an interesting problem that deserves to be investigated.

%

The literature devoted to Hadamard states for Dirac fields on curved spacetimes is much less numerous than the one for scalar fields. 
In particular there seems to be rather few results on  {\em existence} of Hadamard states for Dirac fields.

To our knowledge the first paper proving existence of Hadamard states for Dirac fields in the general case is the recent paper by Murro and Volpe \cite{MV}. It relies on  the construction of an isomorphism between solutions of the Dirac equation  on $(M, g)$ and  on $(M, g_{\rm us})$ respectively,  where $g_{\rm us}$ is an ultrastatic metric having a Cauchy surface in common with $(M, g)$.  This construction is related to the familiar deformation argument of Fulling, Narcowich and Wald \cite{FNW} for bosonic scalar fields. 
 
 Sahlmann and Verch \cite{SV2} have proved that on {\em stationary spacetimes}ground and thermal states are Hadamard, in a general framework which applies in particular to the Dirac case.  Pulling back the ground state for the Dirac operator $D_{\rm us}$ on $(M, g_{\rm us})$ by the above isomorphism, one obtains a state for the Dirac operator $D$ on $(M, g)$, which can be proved to be a Hadamard state, see \cite[Thm. 4.12]{MV}.

Note that the  deformation argument was  already used by d' Antoni and Hollands, see \cite[Sect. V]{DHo} to prove the split property of Hadamard states for  Dirac fields.

Another recent construction  of Hadamard states for Dirac fields  is due to 
Islam and Strohmaier \cite{IS},
where Feynman propagators are studied in detail  for normally hyperbolic operators acting on sections of vector bundles, following the classic approach by Duistermaat and H\"{o}rmander \cite{DH}.

\subsection{Content of this paper}\label{sec0.3}
The goal of this paper is to give an  explicit construction of {\em pure} Hadamard states for Dirac fields on a rather large  class of spacetimes.
Various spacetimes of physical interest like for example the Kerr or Kerr-de Sitter exterior spacetimes, or the Kerr-Kruskal spacetime  describing blackholes fall into the  class of spacetimes that we consider.

We will construct Hadamard states by working on a fixed Cauchy surface $\Sigma$ of $M$, using global {\em pseudodifferential calculus} on $\Sigma$.

To this end we will use the global calculus introduced by Shubin \cite{Sh}, which relies on the notion of {\em bounded geometry}. Let us note that 
the framework  of bounded geometry is often used for global spectral analysis on non compact manifolds, see among many others the papers \cite{AG, AGV}.

In the context of quantum fields on curved spacetimes, Shubin's calculus was first used  in \cite{GOW} to construct pure Hadamard states for Klein-Gordon fields on Lorentzian manifolds of bounded geometry.   

An important fact, probably known to experts, is that on a Lorentzian manifold of bounded geometry possessing a spin structure, 
the spin bundle and hence its associated spinor bundle can, modulo a bundle equivalence, be assumed to be   bundles of bounded geometry.

 This implies that the Dirac operator is itself a differential operator of bounded geometry, and allows to use freely the Shubin's calculus in later steps.

As in \cite{GOW} the first step consists, after fixing a Cauchy surface $\Sigma$  with good properties,  to reduce ourselves using the normal geodesic flow to a product situation where $M= I\times \Sigma$, $I$ some time interval, with the metric $g$ of the form $-dt^{2}+ h(t, \rx)d\rx^{2}$, where $h(t, \rx)d\rx^{2}$ is a time dependent Riemannian metric on $\Sigma$.

The Dirac equation $D\psi=0$ can then be reduced to a time dependent Schroedinger equation
\[
\p_{t}\psi- \i H(t)\psi=0,
\]
where $H(t)$ is a time dependent first order elliptic operator on $\Sigma$. In \cite{GOW} a similar reduction was done by writing the Klein-Gordon equation as a first order system. If $U(t,s)$ denotes the associated Cauchy evolution, one obtains  a pure Hadamard state $\omega$ by constructing a pair of 
selfadjoint projections $P^{\pm}(t)$ with $P^{+}(t)+ P^{-}(t)=\one$, such that
\[
\begin{array}{rl}
i)& U(t, s)P^{\pm}(s)- P^{\pm}(t)U(t,s)\hbox{ is smoothing},\\[2mm]
ii)& \WF(U(\cdot,s)P^{\pm}(s)f)\subset \cN^{\pm}, \ \forall f\in \cE'(\Sigma, S(\Sigma)),
\end{array}
\]
where $S(\Sigma)$ is the restriction of the spinor bundle $S(M)$ to $\Sigma$, $\WF(u)\subset T^{*}M\setminus \zero$ is the {\em wavefront set} of a distribution $u$, and  $\cN^{\pm}$ are the two connected components of the characteristic manifold $\cN$ of the Dirac operator arising in the microlocal definition of Hadamard states.

In \cite{GOW}, the projections $P^{\pm}(t)$ were deduced from a diagonalization of the Cauchy evolution $U(t,s)$ modulo smooth errors, which in turn followed from a factorization of the Klein-Gordon operator. Here the projections are constructed directly by an inductive argument, starting from the spectral projections $\one_{\rr^{\pm}}(H(t))$. Our construction is actually analogous to the construction of {\em adiabatic projections} in linear adiabatic theory, see for example the construction given by Sj\"{o}strand \cite{Sj} in the semiclassical case.

From these projections, one can also obtain an  explicit expression for the {\em Feynman propagator}  associated to the Hadamard state $\omega$.

We conclude our paper by giving two proofs of the existence of Hadamard states for Dirac fields on arbitrary globally hyperbolic spacetimes. The first is by the usual deformation argument of Fulling, Narcowich and Wald \cite{FNW}, the second using a partition of unity on a Cauchy surface is due to \cite{GW1} for Klein-Gordon fields.

\subsection{Plan of the paper}
In Section \ref{sec1} we recall some background on spin structures and Dirac operators on Lorentzian manifolds.  Section \ref{sec1b} is devoted to the quantization of Dirac fields on curved spacetimes and to the definition of Hadamard states. 

In Section \ref{sec2} we recall various definitions related  with  bounded geometry, in particular the notion of Lorentzian manifolds and Cauchy surfaces of bounded geometry, introduced in \cite{GOW}. 

In Section \ref{sec4} we consider the pseudodifferential calculus on manifolds of bounded geometry. We extend several results in \cite{GOW} to pseudodifferential operators acting on sections of bounded vector bundles.

Section \ref{sec5} contains the proof of the main result of the paper, Theorem \ref{maintheorem},  namely the construction of pure Hadamard states for Dirac operators on Lorentzian manifolds of bounded geometry.  Finally in Section \ref{sec7} we give two constructions of Hadamard states for Dirac fields on arbitrary globally hyperbolic spacetimes. The first construction is similar in spirit to the one in \cite{MV}, relying on the deformation argument, but working systematically   with Cauchy data and Cauchy evolutions. 

The second is by a partition of unity argument first used in \cite{GW1} for scalar bosonic fields. 
\subsection{Notations}\label{sec0.-2}
\subsubsection{Lorentzian manifolds}
We use the mostly $+$  signature convention for Lorentzian metrics. All Lorentzian manifolds considered in this paper will be {\em orientable} and connected. 

If $(U_{i})_{i\in \nn}$ is a covering of a manifold $M$, for example associated to a bundle atlas of some bundle $E\xrightarrow{\pi}M$, we set $U_{ij}= U_{i}\cap U_{j}$.

\subsubsection{Bundles}
If  $E\xrightarrow{\pi}M$ is a bundle we denote by $\cinf(M; E)$ resp. $\coinf(M; E)$ the set of smooth resp. smooth and compactly supported sections of $E$.

If $E\xrightarrow{\pi}M$ is a vector bundle of finite rank, we denote by $\cD'(M; E)$ resp. $\cE'(M; E)$ the space of distributional resp. compactly supported distributional sections of $E$.
\subsubsection{Matrices}
Since we will often use frames of vector bundles we will  denote by $\pmb{M}$ a matrix in  $M_{n}(\rr)$ or $M_{N}(\cc)$  and by $M$ the associated endomorphism. 
 \subsubsection{Frames and frame indices} We use the letters $0\leq a \leq d$ for frame indices on $TM$ or $T^{*}M$, and $1\leq a\leq d$ for frame indices on $T\Sigma$ or $T^{*}\Sigma$, if $\Sigma\subset M$ is a space like hypersurface.  If $g$ is a metric on $M$ and $(e_{a})_{0\leq a\leq d}$ is a local frame of $TM$ we set $g_{ab}= e_{a}\dual g e_{b}$ and $g^{ab}= e^{a}\dual g^{-1}e^{b}$, where $(e^{a})_{0\leq a \leq d}$ is the dual frame. 
 
  We use capital letters $1\leq A\leq N$ for frame indices of the spinor bundle $S(M)$.
  
  If $\cF$ is for example a local frame of $TM$  we denote by $\cF\pmb{t}$ the frame obtained by the right action of $\pmb{t}\in M_{n}(\rr)$ on $\cF$.

\subsubsection{Vector spaces}
if $\cX$ is a real or complex vector space, we denote by $\cX'$ its dual. If  $\cX$ is a complex vector space we denote by $\cX^{*}$ its anti-dual, i.e. the space of anti-linear forms on $\cX$ and by  $\bar{\cX}$ its conjugate, i.e. $\cX$ equipped with the complex structure $-\i$. 

A linear map $a\in L(\cX, \cX')$ is a bilinear  form on $\cX$, whose action on pairs of vectors is denoted  by $x_{1}\dual a x_{2}$. Similarly a linear map $a\in L(\cX, \cX^{*})$ is a sesquilinear form on $\cX$, whose action is denoted by $\bar{x}_{1}\dual a x_{2}$. We denote by $a'$, resp. $a^{*}$ the transposed resp. adjoint of $a$.  The space of symmetric resp. Hermitian forms on $\cX$ is denoted by $L_{\rm s}(\cX, \cX')$ resp. $L_{\rm h}(\cX, \cX^{*})$.
\subsubsection{Maps}
We write $f:A\xrightarrow{\sim}B$ if $f: A\to B$ is a bijection. We use the same notation if $A, B$ are topological spaces resp. smooth manifolds, replacing bijection by homeomorphism, resp. diffeomorphism.
 \section{Spin structures and Dirac operators on Lorentzian spacetimes}\init\label{sec1}
 We recall some definitions about spin structures and Dirac operators, see \cite[chap. I , chap.  II sects 1,3,4,  5]{LM}. Another nice exposition is found in Trautman \cite{T}.
 \subsection{Background}\label{sec1.1}
 We denote by  $\rr^{1, d}$  the Minkowski spacetime, i.e. $\rr^{1+d}$ equipped with the bilinear form $x\dual \eta x= - t^{2}+ \rx^{2}$, $x= (t, \rx)$.
 We will always assume that $n=1+d$ is {\em even}.    We denote by $\SO^{\uparrow}(1, d)\subset \SO(1, d)\subset {\rm O}(1, d)$ the restricted Lorentz group. 
  \subsubsection{Clifford algebras}\label{sec1.1.1}
 The {\em Clifford algebra} $\Cl(1, d)$ is the real algebra  generated by elements $\one, \gamma(x),x\in \rr^{1, d}$ with relations
 \beq\label{cliff-conv}
 \gamma(x)\gamma(x')+ \gamma(x')\gamma(x)= 2 x\dual \eta x'\one.
 \eeq
 For $r\in {\rm O}(1, d)$ one sets  $\hat{r}(\gamma(x))= \gamma(rx)$ $x\in \rr^{1, d}$.  The map $\hat{r}$ extends uniquely as an automorphism of $\Cl(1, d)$. We set
 \beq\label{e1.-2}
\alpha:  {\rm O}(1, d)\ni r\mapsto \hat{r}\in Aut(\Cl(1, d))
 \eeq which is  a group morphism.
 \subsubsection{Spin groups}\label{sec1.1.2}
The {\em Spin group} $\Spin(1, d)\subset \Cl(1, d)$ is
 \[
 \Spin(1, d)= \{\gamma(x_{1})\cdots \gamma(x_{2p}): x_{i}\dual \eta x_{i}= \pm 1, p\in\nn\}
 \]
 and the {\em restricted Spin group} $\Spin^{\uparrow}(1, d)$ is the connected component of $\one$ in $\Spin(1, d)$. One can show that $\gamma(x_{1})\cdots \gamma(x_{2p})$ belongs to $\Spin^{\uparrow}(1, d)$ iff the number of indices $i$ with $x_{i}\dual \eta x_{i}= -1$ is even.

The map  $Ad: \Spin^{\uparrow}(1, d)\to \SO^{\uparrow}(1, d)$ defined by
 \[
 a \gamma(x)a^{-1}\eqdef \gamma(Ad(a)x), \ a\in \Spin^{\uparrow}(1, d)
 \]
 is a two-sheeted covering  and one has $Ad^{-1}(\{r\})= \{a, -a\}$ for $Ad(a)= r$.
 \subsubsection{Representations}\label{sec1.1.3}
 Since $n$ is even, there exist a (unique up to isomorphism) faithful and irreducible representation 
 \beq\label{e1.-5}\rho_{0}: \Cl(1, d)\to M_{N}(\cc), \ N= 2^{\frac{n}{2}}.
 \eeq
We equip $\rr^{1, d}$ and $\cc^{N}$ with their canonical bases $(u_{a})_{0\leq a\leq d}$, $(v_{A})_{1\leq A\leq N}$
and tacitly identify $\SO^{\uparrow}(1, d)$, $\Cl(1, d)$ and $\Spin^{\uparrow}(1, d)$ with their images in $M_{n}(\rr)$ or $M_{N}(\cc)$.

We denote hence by $\pmb{\gamma}_{a}\in M_{N}(\cc)$ the matrix of $\rho_{0}(\gamma(u_{a}))$  and set
\[
\pmb{\gamma}(v)= \pmb{\gamma}_{a}v^{a}, \ v= v^{a}e_{a	}\in \rr^{1, d}.
\]
There exists a Hermitian matrix  $\pmb{\beta}\in M_{N}(\cc)$ such that
\begin{equation}
\label{e2.2}
\begin{array}{l}
\pmb{\gamma}(v)^{*} \pmb{\beta}= -\pmb{\beta}\pmb{\gamma}(v),\ v\in \rr^{1, d},\\[2mm]
 \i \pmb{\beta}\pmb{\gamma}(v)>0\hbox{ if }v\hbox{ is  time-like future directed}. 
\end{array}
\end{equation}
If $n\in \{2, 4\}$ mod $8$ there exists  a real matrix $\pmb{\kappa}\in M_{N}(\cc)$  such that
\[
\pmb{\kappa}^{2}= \one, \ \pmb{\kappa}\pmb{\gamma}(v)= \pmb{\gamma}(v)\pmb{\kappa}, \ v\in \rr^{1, d}.
\]
One can show that $\Spin^{\uparrow}(1, d)$ is the set of $\pmb{a}\in M_{N}(\cc)$ such that
\begin{equation}
\label{e2.3}
\begin{array}{rl}
(i)&\pmb{a}^{*}\pmb{\beta}\pmb{a}=\pmb{\beta}, \ \pmb{a} \pmb{\kappa}=  \pmb{\kappa}\overline{\pmb{a}},\\[2mm]
(ii)&\pmb{a} \pmb{\gamma}(v)\pmb{a}^{-1}=\pmb{\gamma}(Ad(\pmb{a})v), \ v\in \rr^{1, d}.
\end{array}
\end{equation}
\subsection{Embedding of $\rr^{d}$ into $\rr^{1, d}$}\label{sec1.1b}
Let $\rr^{d}$ be the Euclidean space, equipped with the bilinear form $\rx\dual \delta \rx= \rx^{2}$.  We denote by $\Cl(d)$, $\SO(d)$, $\Spin(d)$ the analogous Clifford algebra, special orthogonal group and spin group with $\eta$ replaced by $\delta$. 

We consider the  isometric embedding ${\rm i}: \rr^{d}\ni \rx\mapsto (0, \rx)\in \rr^{1, d}$. Note that it is different from the one considered in \cite[Sect. 2]{BGM} since  the range of ${\rm i}$ is space-like in our case. It induces a morphism  from $\Cl(d)$ to $\Cl(1, d)$, still denoted by ${\rm i}$ for simplicity of notation. It is easy to see that this morphim is injective.   

In fact setting  $\pmb{E}_{I}= \prod_{i\in I}\pmb{\gamma}_{i}$ for $I\subset \{0, \dots, d\}$, then $\{\pmb{E}_{I}\}_{I\subset \{0, \dots, d\}}$ is a basis of $\Cl(1, d)$ while $\{\pmb{E}_{I}\}_{I\subset \{1, \dots, d\}}$ is a basis of $\Cl(d)$ and ${\rm i}\pmb{E}_{I}= \pmb{E}_{I}$ for $I\subset \{1, \dots, d\}$ which proves the injectivity of ${\rm i}$.
 
 Using the characterization of $\Spin^{\uparrow}(1, d)$ in \ref{sec1.1.2} and the same characterization of $\Spin(d)$ we obtain that ${\rm i}: \Spin(d)\to \Spin^{\uparrow}(1, d)$. 
 
 Note that we have also an injective morphism
  \beq\label{tru}
  \tilde{\rm i}: \SO(d)\ni \pmb{o}\mapsto \mat{1}{0}{0}{\pmb{o}}\in \SO^{\uparrow}(1, d)
  \eeq 
  with a commutative diagram 
 \[
  \begin{tikzcd}
\SO(d)\arrow[r, "\tilde{\rm i}"] \arrow[d, "\alpha"] &\SO^{\uparrow}(1, d) \arrow[d, "\alpha"]    \\
 \Cl(d) \arrow[r, "\rm{i}"]            & \Cl(1, d)              
\end{tikzcd}
\]
 where $\alpha$ denotes the action of $\SO(d)$, resp. $\SO(1, d)$ on $\Cl(d)$, resp. $\Cl(1, d)$.

\begin{lemma}\label{lemma1.1}
Let $\pmb{s}\in \Spin^{\uparrow}(1, d)$ with $Ad(\pmb{s})\in \tilde{\rm i}(\SO(d))$. Then there  exists a unique $\pmb{\tilde{s}}\in \Spin(d)$ such that $\pmb{s}= \rm{i}(\pmb{\tilde{s}})$.
\end{lemma}
\proof
We write $\pmb{s}= \sum_{I\subset \{1, \dots, d\}}(\lambda_{I}\pmb{\gamma}_{0}+ \mu_{I})\pmb{E}_{I}$. Since $Ad(\pmb{s})e_{0}= e_{0}$  we have $\pmb{s}\pmb{\gamma}_{0}= \pmb{\gamma}_{0}\pmb{s}$, which using the Clifford relations implies that $\lambda_{I}=0$ for all $I$ and hence $\pmb{s}\in \rm{i}(\Spin(d))$. \qed

 \subsection{Lorentzian manifolds}\label{sec1.2}
 
 \subsubsection{Spacetimes}\label{sec1.2.1}
 We recall that a  spacetime $(M , g)$ is an orientable and time-orientable Lorentzian manifold.  We will always assume that $M$ is connected. 
 \subsubsection{The bundle $P\SO^{\uparrow}(M, g)$}\label{sec1.2.2}
 Let $M$ be a spacetime.  It is well-known that $TM$ admits local oriented and time oriented orthonormal frames. For completeness let us sketch the proof of this fact:
 
 let  $(\mathcal{F}_{i})_{i\in \nn}$  a family of local frames of $TM$ over $(U_{i})_{i\in \nn}$, with $\mathcal{F}_{i}= (f_{i, a})_{0\leq a\leq d}$, and 
$\pmb{g}_{i, ab}= f_{i, a}\dual g f_{i, b}$.  Without loss of generality we can assume that $\mathcal{F}_{i}$ is direct and $f_{i, 0}$ is future directed.  We choose a neighborhood $V_{0}$ of $\pmb{\eta}$ as in Lemma \ref{lemma-app.1} such that $\mathcal{F}_{i}\pmb{t}$  is direct and $\pmb{t}f_{i, 0}$ is future directed for all $\pmb{t}\in U_{0}= F(V_{0})$, where the map $F$ is constructed in Lemma \ref{lemma-app.1}.

We set then $\pmb{t}_{i}\defeq  F(\pmb{g}_{i})$, $\mathcal{E}_{i}\defeq  \mathcal{F}_{i}\pmb{t}_{i}$ and $(\mathcal{E}_{i})_{i\in \nn}$ is a family of local orthonormal oriented and time oriented frames of $TM$.

 We denote by $P\SO^{\uparrow}(M, g)$ the $\SO^{\uparrow}(1, d)$-principal bundle over $M$ of local oriented and time oriented orthonormal frames of $TM$.

Equivalently one can define   $P\SO^{\uparrow}(M, g)$ as the $\SO^{\uparrow}(1, d)$-principal bundle over $M$ with transition functions 
 \beq\label{defdeoij}
 \pmb{o}_{ij}= \pmb{t}_{j}^{-1}\circ \pmb{t}_{i}: U_{ij}\to \SO^{\uparrow}(1, d).
 \eeq
 \subsubsection{The bundle $\Cl(M, g)$}\label{sec1.2.3}
 From the bundle $P\SO^{\uparrow}(M, g)$ and the map introduced in \eqref{e1.-2} we obtain by the associated bundle construction the {\em Clifford bundle} $\Cl(M, g) = P\SO^{\uparrow}(M, g)\times_{\alpha}\Cl(1,d)$, which is a bundle of algebras with typical fiber $\Cl(1, d)$. Equivalently one can define $\Cl(M, g)$ by the transition maps
 \[
 \hat{\pmb{o}}_{ij}: U_{ij}\to Aut(\Cl(1, d)),
 \]
 for $\pmb{o}_{ij}$ defined in \eqref{defdeoij}.

\subsection{Spin structures}\label{sec1.3}
 \begin{definition}
Let $(M, g)$ be a spacetime. A {\em spin structure} on $(M, g)$ is a $\Spin^{\uparrow}(1, d)$-principal bundle over $M$ denoted   $P\Spin(M, g)$ such that there exists a morphism of principal bundles $\chi: P\Spin(M, g)\to P\SO^{\uparrow}(M, g)$ such that the following diagram commutes:\beq\label{e1.-3}
 \begin{tikzcd}
\Spin^{\uparrow}(1,d) \arrow[r] \arrow[dd, "Ad"] &P\Spin(M, g) \arrow[rd, "\pi'"] \arrow[dd, "\chi"] &   \\
                   &                         & M. \\
\SO^{\uparrow}(1, d) \arrow[r]            & P\SO^{\uparrow}(M, g) \arrow[ru, "\pi"]            &  
\end{tikzcd}
\eeq
\end{definition}
\subsubsection{Existence and uniqueness of spin structures}\label{sec1.3.0}
For the reader's convenience, let us recall well-known results on the existence and uniqueness of spin structures.

A Lorentzian manifold admits  a spin structure if and only if  its  second Stiefel--Whitney class $w_{2}(TM)$ is trivial, see \cite{Mi, Na}.  If $\dim M=4$ this is also equivalent to the fact that $M$ is parallelizable,  see \cite{spinors1, spinors2}. 
It admits a {\em unique} spin structure if in addition  its  first Stiefel--Whitney class $w_{1}(M)$ is trivial, which is equivalent to the fact that $M$ is orientable, see e.g.~\cite{Na}.

In our situation,   $M$ is orientable hence spin structures on $(M, g)$ are unique  if they exist. If $(M, g)$ is globally hyperbolic, then $M$ is diffeomorphic to $\rr\times \Sigma$, where $\Sigma$ is a smooth Cauchy surface. The orientation and time orientation of $M$  induce an orientation of $\Sigma$, hence $\Sigma$ is orientable.

Since any orientable $3$-manifold is parallelizable, this implies that if $(M, g)$ is globally hyperbolic and  $n=4$,  $M$ is parallelizable and by the above facts it admits a  {\em unique} spin structure.

\subsection{Spinor bundle}\label{sec1.3.1}
From the bundle $P\Spin(M, g)$ and the map \eqref{e1.-5} we obtain by the associated bundle construction the {\em spinor bundle} $S(M) = P\Spin(M,g)\times_{\rho_0}\cc^N$ which is a vector bundle with typical fiber $\cc^{N}$ and   the same transition maps 
\[
\pmb{s}_{ij}: U_{ij}\to M_{N}(\cc)
\]
 as $P\Spin(M, g)$.

The spinor bundle $S(M)$ inherits a lot of extra structures which we now recall.  

\subsubsection{Action of $\Cl(M, g)$ on $S(M)$}\label{sec1.3.2}
We can define a morphism of bundles of algebras $\rho : \Cl(M,g)\to L(S(M))$ as follows
\[\rho([\chi(P),a])\cdot [P,v]=[P,\rho_0(a)\cdot v],\]
where $P \in P\Spin(M,g)$, $a\in \Cl(1,d)$, and $v\in \cc^N$.

We could also define it directly from transition maps.
From \eqref{e2.3} (ii) we obtain that if $\pmb{s}_{ij}$ are the transition maps of $S(M)$ and  $Ad(\pmb{s}_{ij})\eqdef \pmb{o}_{ij}$,  then 
\[
\pmb{s}_{ij}\pmb{\gamma}(v)\pmb{s}_{ij}^{-1}= \pmb{\gamma}(\pmb{o}_{ij}v)= \alpha(\pmb{o}_{ij})\pmb{\gamma}(v), \ v\in \rr^{1, d}.
\]
Therefore the morphism $\rho_{0}: \Cl(1, d)\to M_{N}(\cc)$ induces a 
morphism of bundles of algebras:
\[
\rho: \Cl(M, g)\to L(S(M))
\]
For $v\in \cinf(M; TM)$ we set
\begin{equation}
\label{e2.4}
\gamma(v)\defeq  \rho(v), \hbox{ where }v\in \cinf(M; TM)\subset \cinf(M; \Cl(M, g)).
\end{equation}
The map $\gamma(v)$ is usually called the Clifford multiplication by $v$ and often denoted by $v\cdot$.
\subsubsection{Time positive Hermitian structure}\label{sec1.3.3}
We can define a Hermitian structure $\beta \in\cinf(M; L_{\rm h}(S(M), S(M)^{*}))$ by
\[[P,u]\cdot \beta [P,v] = u\cdot \pmb{\beta}v,\]
where $P \in P\Spin(M,g)$, $u,v\in \cc^N$ and $\pmb{\beta}$ is in \eqref{e2.2}.

Equivalently using  \eqref{e2.3} and  denoting by $T_{i}: S(M)\cap \pi^{-1}(U_{i})\to U_{i}\times \cc^{N}$ the local trivializations of $S(M)$, we can define  $\beta$ by:
\[
\beta= T_{i}^{*}\pmb{\beta}T_{i}\hbox{ over }U_{i}.
\]
\subsubsection{Charge conjugation}\label{sec1.3.4}
 Similarly we can define a charge conjugation $\kappa$ belonging to $\cinf(M; L(S(M), \bar{S(M)}))$ by
 \[\kappa [P,v] = [P,\pmb{\kappa} v]\]
where $P \in P\Spin(M,g)$, $v\in \cc^N$ and $\pmb{\kappa}$ is in \eqref{e2.3}.

Equivalently using  \eqref{e2.3}  we can define $\kappa$ by\[
 \kappa= T_{i}^{-1}\pmb{\kappa}T_{i}\hbox{ over }U_{i},
\]
where above we still denote by $\pmb{\kappa}$ the anti-linear map $\cc^{N}\ni z\mapsto \pmb{\kappa}\bar{z}$.
\subsubsection{Properties}\label{sec1.3.5}
Using \eqref{e2.2}, \eqref{e2.3} we obtain
\begin{equation}
\label{e2.5}
\begin{array}{l}
\gamma(v)^{*} \beta= - \beta\gamma(v), \ v\in\cinf(M; TM),\\[2mm]
 \kappa\gamma(v)= \gamma(v)\kappa, \ v\in \cinf(M; TM)\\[2mm]
 \i \beta\gamma(e)>0\hbox{ if }e\in \cinf(M; TM)\hbox{ is future directed}. 
\end{array}
\end{equation}
\subsubsection{Local frame of $S(M)$ associated to a local frame of $TM$}\label{sec1.3.6}
Let   $U\subset M$ a chart open set and $\mathcal{F}= (e_{a})_{0\leq a\leq d}$   an orthonormal oriented and time oriented local frame   of $TM$ over $U$, i.e. a local section of $P\SO^{\uparrow}(M, g)$.  If $S$ is a local section of $P\Spin^{\uparrow}(M, g)$ such that $Ad(S)= \cF$, then using  the concrete representation $\rho_{0}$ and the canonical basis $(v_{A})_{1\leq A\leq N}$  of $\cc^{N}$ we obtain a local frame $\mathcal{B}= (E^{A})_{1\leq A\leq N}$ of $S(M)$ over $U$.

The  spin frame $\mathcal{B}$ is said {\em associated} to the vector frame $\mathcal{F}$.

Modulo a bundle isomorphism of $S(M)$ over $U$ we can assume that the matrices of  $\beta$, resp. $\kappa$, $\gamma(e_{a})$ in the frame $ (E^{A})_{1\leq A\leq N}$ are given by $\pmb{\beta}$ resp.  $\pmb{\kappa}$, $\pmb{\gamma}_{a}$.

\subsubsection{Spin connection}\label{sec1.3.7}
The Levi-Civita  connection on $P\SO^{\uparrow}(M, g)$ lifts through $\chi$ in (\ref{e1.-3}) to a connection on $P\Spin(M, g)$  and hence on $S(M)$ called the {\em spin connection}. 

Denoting by $\nabla$, resp. $\nabla^{S}$ the associated covariant derivatives, we have:
\beq\label{e10.4}
\begin{array}{rl}
i)&\nabla_{X}^{\cS}(\gamma(Y)\psi)= \gamma(\nabla_{X}Y)\psi+ \gamma(Y)\nabla_{X}^{\cS}\psi,\\[2mm]
ii)&X(\bar{\psi}\dual \beta \psi)= \bar{\nabla_{X}^{\cS}\psi}\dual \beta \psi+ \bar{\psi}\dual \beta \nabla_{X}^{\cS}\psi,\\[2mm]
iii)&\kappa \nabla_{X}^{\cS}\psi= \nabla_{X}^{\cS}\kappa \psi,
\end{array}
\eeq
for all $X, Y\in \cinf(M; TM)$ and $\psi\in \cinf(M; \cS(M))$.

  A concrete expression for the spin covariant derivative $\nabla^{S}$ is as follows: 
  
  let  $(e_{a})_{0\leq a\leq d}$  an orthonormal oriented and time oriented local frame of $TM$ over $U$ and $(e^{a})_{0\leq a\leq d}$ the dual frame. We   set 
 \[
\nabla_{b}\defeq  \nabla_{e_{b}}, \ \Gamma^{c}_{ba}\defeq  \nabla_{b}e_{a}\dual  e^{c}\]
hence
\[
\nabla_{u}v= u^{b}(\p_{b}v^{a}+ \Gamma^{a}_{bc}v^{c})e_{a}.
\]
where we denote
\[
 u^{a}\defeq  u\dual e^{a},  \ \p_{b}f\defeq  e_{b}\dual df
\hbox{ for }u\in \cinf(U, TM),  f\in \cinf(U).
\]

If $ (E_{A})_{1\leq A\leq N}$ is the associated spin frame  and $(E^{A})_{1\leq A\leq N}$ the dual frame, we set 
\[
 \psi^{A}\defeq  \psi\dual E^{A}, \ M^{A}_{B}\defeq  E^{A}\dual M E_{B} \hbox{ for }\psi\in \cinf(U; S(M)); \ M\in \cinf(U, L(S(M))).
\]
Then we have:
\beq\label{e2.7}
\nabla^{S}_{u}\psi\defeq  u^{b}(\p_{b}\psi^{A}+\sigma^{A}_{bC}\psi^{C})E_{A},
\eeq
where
\[
\sigma_{b}\defeq  \frac{1}{4}\Gamma^{a}_{bc}\gamma_{a}\gamma^{c},\ \gamma_{a}= \gamma(e_{a}), \ \gamma^{b}= \eta^{ab}\gamma_{a}.
\]
\subsection{Spin structures on cartesian products}\label{spinprod}
Let $M= I_{t}\times \Sigma_{\rx}$ where $I\subset \rr$ is an open interval and $\Sigma$ an orientable $d$-dimensional manifold.  We fix  the Lorentzian metric  on $M$:
\[
g= -dt^{2}+ h(t,\rx)d\rx^{2},
\]
where $I\ni t\mapsto h(t, \rx)d\rx^{2}$ is a family of Riemannian metrics on $\Sigma$,   and  equip $M$  with the natural orientation and time orientation obtained from the orientation of $\Sigma$.  Let us set $h_{t}(\rx)d\rx^{2}= h(t, \rx)d\rx^{2}$.
 
Let $(V_{i})_{i\in \nn}$ an open covering of $\Sigma$, $\cF_{i}= (e_{i, a})_{1\leq a\leq d}$ local oriented orthonormal frames for $h_{0}$ over $V_{i}$
and $\pmb{o}_{ij}: V_{ij}\to \SO(d)$ the associated transition functions. Let $\cF_{i}(t)= (e_{i, a}(t))_{1\leq a\leq d}$ the oriented orthonormal frames for $h_{t}$ over $V_{i}$ obtained by parallel transport of $\cF_{i}$ with respect to $\p_{t}$ for the metric $g$, and $\pmb{o}_{ij}(t): V_{ij}\to \SO(d)$ the associated transition functions.  Since  the entries of $\pmb{o}_{ij}(t)$ are given by $e_{i, a}(t)\dual h_{t}e_{j, b}(t)= e_{i, a}(t)\dual g e_{j, b}(t)$ we obtain that $\pmb{o}_{ij}(t)=\pmb{o}_{ij}$ is independent on $t$. 

Note also that $\tilde{i}(\pmb{o}_{ij}): \rr\times V_{ij}\to \SO^{\uparrow}(1, d)$ are the transition functions of $P\SO^{\uparrow}(M, g)$, since if we set  
$e_{0}= \p_{t}$, then $\cE_{i}= (e_{i, a}(\cdot))_{0\leq a\leq d}$ are local oriented and time oriented orthonormal frames for $g$ over $U_{i}= I\times V_{i}$.
\subsubsection{Spin structures}

We use the notation in Subsect. \ref{sec1.1b}.
Assume that $(M, g)$ has a spin structure, which is then unique since $M$ is orientable.  By Prop. \ref{propadd1.1} we can assume that the associated transition maps $\pmb{s}_{ij}: U_{ij}\to \Spin\uparrow(1, d)$ satisfy $Ad(\pmb{s}_{ij})= \tilde{\i}(\pmb{o}_{ij})$, and are hence independent on $t$.

By Lemma \ref{lemma1.1} there exists unique $\pmb{\tilde{s}}_{ij}: U_{ij}\to \Spin(d)$, independent on $t$,  such that  $\pmb{s}_{ij}= {\rm i} (\pmb{\tilde{s}}_{ij})$ and hence   we obtain a spin structure on $(\Sigma, h_{t})$ for $t\in I$.

Conversely if $(\Sigma, h_{0})$ has a spin structure with transition functions $\pmb{\tilde{s}}_{ij}: V_{ij}\to \Spin(d)$, the  transition functions \[
\pmb{s}_{ij}= {\rm i}(\tilde{\pmb{s}}_{ij}): U_{ij}\to \Spin^{\uparrow}(1, d),
\]
 define  a spin structure on $(M, g)$, since $Ad(\pmb{s}_{ij})= \tilde{i}(\pmb{o}_{ij})$. 
\subsubsection{Spinor bundles}
Let  $S(M)$ the  spinor bundle on $(M, g)$ and $S_{t}(\Sigma)$ its restriction to $\{t\}\times \Sigma$. Since the transition maps $\pmb{s}_{ij}$ are independent on $t$, $S_{t}(\Sigma)$ is independent on $t$ and will hence be denoted by $S(\Sigma)$. 
\subsection{Dirac operators}\label{sec1.4}
\begin{definition}\label{def1.2}
Let $(M, g)$ a spacetime with a spin structure and $S(M)$ the associated spinor bundle.  A {\em Dirac operator} on $M$ is a  differential operator $D$ acting on $\cinf(M; S(M))$  such that if  $(e_{a})_{0\leq a\leq d}$ is a local frame over $U\subset M$ then
\[
D= g^{ab}\gamma(e_{a})\nabla^{S}_{e_{b}}+ m \hbox{ over }U
\]
where  $m\in \cinf(M, L(S(M)))$ is a section such that $m^{*}\beta = \beta m$. 
\end{definition}
We  set 
\[
\slashed{D}\defeq  g^{ab}\gamma(e_{a})\nabla^{S}_{e_{b}}.
\]

\subsubsection{Selfadjointess of Dirac operator}
Let $\psi_1,\psi_2 \in \cinf(M; S(M))$. 
One defines the $1$-form $J(\psi_1,\psi_2) \in \cinf(M;T^*M)$ by
\[J(\psi_1,\psi_2)\cdot X \defeq  \overline{\psi_1} \cdot \beta \gamma(X) \psi_2, X\in \cinf(M;TM),
\]
and checks using \eqref{e10.4} that 
\[
\nabla^{\mu} J_{\mu}(\psi_1,\psi_2) = - \overline{D\psi_1}\cdot \beta \psi_2 + \overline{\psi_1}\cdot \beta D \psi_2,\ \psi_i \in \cinf(M; S(M)).
\]	
This implies easily the following proposition.
\begin{proposition} \label{prop1.0} The
	Dirac operator $D$ is formally self-adjoint on $\coinf(M; S(M))$ with respect to the Hermitian form:
	\beq\label{burk}
	(\psi_1 | \psi_2)_M \defeq  \int_M \psi_1 \cdot \beta \psi_2 dVol_g.\eeq
\end{proposition}

\subsubsection{Conformal transformations}\label{sec1.4.1}
We recall the well-known covariance of Dirac equations under conformal transformations   see \cite[Lemma 5.27]{LM} or \cite{Hi}.
If $\tilde{g}= e^{2u}g$ for $u\in \cinf(M)$, we then have the natural isomorphism
\[
G_u:\ \begin{array}{l}
P\SO^{\uparrow}(M,g) \to P\SO^{\uparrow}(M,\tilde{g})\\[2mm]
\cF= (e_{a})_{0\leq a\leq d}\mapsto (\e^{-u}e_{a})_{0\leq a\leq d}= \e^{- u} \cF.
\end{array}
\] 
$G_{u}$ induces an isomorphism still denoted by $G_{u}$ between $\Cl(M, g)$ and $\Cl(M, \tilde{g})$:
\[
G_{u}:\ \begin{array}{l}
\Cl(M,g) \to \Cl(M,\tilde{g})\\[2mm]
[\cF,  a]\mapsto [\e^{-u}\cF, a].
\end{array}
\]
We obtain  a spin structure on $(M,\tilde{g})$ from the one on $(M,g)$ by
\[
\begin{tikzcd}
&\arrow[ldd, "\chi"] P\Spin(M, g) \arrow[dd, "\tilde{\chi}\defeq G_u \circ \chi"] &   \\
 & \\
P\SO^{\uparrow}(M, g)  \arrow[r,"G_u"]  & P\SO^{\uparrow}(M, \tilde{g}) 
\end{tikzcd}
\]
which leads to the same spinor bundles $S(M)$.  Let  $\tilde{\gamma}$, $\tilde{\beta}$, $\tilde{\kappa}$ be as in Subsect. \ref{sec1.3.1} with $g$ replaced by $\tilde{g}$. Clearly we have
\[
\tilde{\beta}= \beta, \ \tilde{\kappa}= \kappa.
\]
We define $\tilde{\rho}: \Cl(M,\tilde{g})\to End(S(M))$ such that the following diagram is commutative
\[
\begin{tikzcd}
P\SO^{\uparrow}(M, g)\times \cc^N \arrow[rdd, "\rho"]    \arrow[r,"G_u"]  & P\SO^{\uparrow}(M, \tilde{g})\times \cc^N \arrow[dd, "\tilde{\rho}"]\\
& \\
 &End(S(M)).
\end{tikzcd}
\]
Considering $TM$ as the associated bundle $TM= P\SO(M, g)\times\rr^{n}$ resp. $TM= P\SO(M, \tilde{g})\times\rr^{n}$ leads to the map $TM\ni X\mapsto \e^{-u}X\in TM$. This leads to $\tilde{\gamma}(\e^{- u}X)= \gamma(X)$ hence:
\[
\tilde{\gamma}(X)= \e^{u}\gamma(X), \ X\in TM.
\]
If $\widetilde{\nabla}^{S}$ is the spinorial covariant derivative for $\tilde{g}$ we have see e.g. \cite[Prop. 4.2.1]{Hi}:
\[\widetilde{\nabla}_X \psi = \nabla_X \psi + \frac{1}{2}\gamma(X)\gamma({\rm grad}(u))\cdot\psi - \frac{1}{2}(X\cdot du) \psi.
\]
(The change of sign in the second term in the rhs in comparison with \cite[Prop. 4.2.1]{Hi} comes from our convention for $\Cl(1, d)$).

 Concerning the Dirac operators we have:
\begin{equation}
\label{econfo.1}
 \tilde{\slashed{D}}= e^{-\frac{n+1}{2}u}\slashed{D}e^{\frac{n-1}{2}u}.
\end{equation}

If we introduce the map
\[
W : \psi\in\cinf_0(M;S(M))\mapsto e^{\frac{n-1}{2}u}\psi,
\]
we have
\[
(\psi_1|W\psi_2)_{(M,g)} = (W^* \psi_1 | \psi_2)_{(M,\tilde{g})}, \ W^*\psi = e^{-\frac{n+1}{2}u}\psi,
\]
where we denote by $(\cdot| \cdot)_{(M, g)}$ resp. $(\cdot | \cdot)_{(M, \tilde{g})}$ the scalar product in \eqref{burk} for $g, \beta$ resp. $\tilde{g},\tilde{\beta}$.
It follows that  (\ref{econfo.1}) can be rewritten as
\begin{equation}
\label{econfo.2}
W^*DW = \tilde{\slashed{D}}+e^{-u}m.
\end{equation}
\section{Quantization of the Dirac equation on curved spacetimes}\label{sec1b}\init
In this section we recall the algebraic quantization of the Dirac equation on curved spacetimes, due to Dimock \cite{Di}.  It is useful to start with a rather general framework see e.g. \cite[Sects 12.5 and 17.2]{DG}).
\subsection{$\CAR*$-algebras and quasi-free states}\label{sec1b.1}
\begin{definition} \label{def1.3}
Let $(\cY,\nu)$ be a (complex) pre-Hilbert space. The $\CAR$ $*${\em -algebra} over $(\cY,\nu)$, denoted by $\CARY$, is the unital complex $*$-algebra generated by elements $\psi(y)$,$\psi^*(y)$, $y \in \cY$ with the relations
\begin{align*}
	\psi(y_1 + \lambda y_2) &= \psi(y_1) + \overline{\lambda}\psi(y_2),\\
\psi(y_1 + \lambda y_2) &= \psi(y_1) + \overline{\lambda}\psi(y_2), y_1,y_2\in \cY, \lambda \in\cc\\
\psi^*(y_1 + \lambda y_2) &= \psi^*(y_1) + \lambda\psi^*(y_2),\\
[\psi(y_1),\psi(y_2)]_+ &= [\psi^*(y_1),\psi^*(y_2)]_+ = 0,\\
[\psi(y_1),\psi^*(y_2)]_+ &= \overline{y_1}\cdot \nu y_2 \id, \ y_1,y_2\in \cY\\
\psi(y)^* &= \psi^*(y), \ y\in \cY,
\end{align*}
where $[A,B]_+ = AB+BA$ is the anti-commutator.
\end{definition}
\subsubsection{Pure quasi-free states}
\begin{definition}
A {\em state} on a $*$-algebra $\fA$ is a linear map $\omega : \fA \rightarrow \cc$ which is normalized, i.e. $\omega(\id) = 1$, and positive, i.e. $\omega(A^*A) \geq 0$ for $A \in \fA$.
The set of states on $\fA$ is a convex set. Its extreme points are called {\em pure states}.
\end{definition}

\begin{definition}
A state $\omega$ on $\CARY$ is a gauge invariant {\em quasi-free state} if

\begin{align*}
\omega(\prod_{i=1}^{n} \psi^*(y_i) \prod_{j=1}^{m} \psi(y_j)) &=0, \text{for $n \not = m$} \\
\omega(\prod_{i=1}^{n} \psi^*(y_i) \prod_{j=1}^{n} \psi(y_j)) &= \sum_{\sigma \in \mathcal{S}_n} sgn(\sigma) \prod_{i=1}^n \omega( \psi^*(y_i)\psi(y_{\sigma(i)}))\end{align*}
\end{definition}
A quasi-free state is characterized by its covariances $\lambda^{\pm} \in L_{\rm h}(\cY,\cY^*)$ defined by
\[\omega(\psi(y_1)\psi^*(y_2))\eqdef \overline{y_1}\cdot \lambda^+ y_2,\ \ \omega(\psi^*(y_1)\psi(y_2))\eqdef \overline{y_1}\cdot \lambda^- y_2, \ \ y_1,y_2 \in \cY\]
The following two results are well-known, see e.g. \cite[Sect. 17.2.2]{DG}.
\begin{proposition}
Let $\lambda^{\pm} \in L_{\rm h}(\cY,\cY^*)$. Then the following statements are equivalent :
\begin{enumerate}
\item $\lambda^{\pm}$ are the covariances of a gauge invariant quasi-free state on $\CARY$,
\item $\lambda^{\pm} \geq 0$ and $\lambda^+ + \lambda^- = \nu$.
\end{enumerate}
\end{proposition}

\begin{proposition}\label{hurlub}
A quasi-free state $\omega$ on $\CARY$ is pure if and only if there exist projections $c^{\pm} \in L(\cY)$ such that
\[\lambda^{\pm} = \nu \circ c^{\pm}, c^++c^-= \id\]
\end{proposition}

\subsection{Quantization of the Dirac equation}\label{sec1b.2}
Assume now that $(M,g)$ is a globally hyperbolic spacetime with  a spin structure. 
We denote by $\Sol(D)$ the space of smooth, space compact solutions of the Dirac equation
\[
D\psi=0.
\]

\subsubsection{Retarded/advances inverses}
Since $(M,g)$ is globally hyperbolic, $D$ admits  unique retarded/advanced inverses $G_{\rm ret/adv} :\coinf(M; S(M)) \rightarrow \cinf(M; S(M))$  see e.g. \cite{Di, Mu} such that

\begin{equation*}
  \left\{
      \begin{aligned}
       & D G_{\rm ret/adv} = G_{\rm ret/adv} D = \id, \\
       & \text{supp } G_{\rm ret/adv} u \subset J_{\pm}(\text{supp }u), u\in\coinf(M; S(M)),
      \end{aligned}
    \right.
\end{equation*}
where $J_{\pm}(K)$ is the future/past causal shadow of $K\subset M$.

Using Prop.  \ref{prop1.0} and the uniqueness of $G_{\rm ret/adv}$ we obtain that
\[G_{\rm ret/adv}^*=G_{\rm adv/ret}\]
where the adjoint is computed with respect to $(\cdot|\cdot)_M$.
Therefore, the causal propagator 
\[G \defeq  G_{\rm ret} - G_{\rm adv}\]
satisfies

\begin{equation*}
  \left\{
      \begin{aligned}
       & D G = G D = 0, \\
       & \text{supp } G u \subset J(\text{supp }u), u\in\coinf(M; S(M)), \\
       & G^*=-G,
      \end{aligned}
    \right.
\end{equation*}
where $J(K)= J_{+}(K)\cup J_{-}(K)$ is the causal shadow of $K\subset M$.
\subsubsection{Cauchy problem}

Let $\Sigma \subset M$ be a smooth, space-like Cauchy surface,  $n$ its future directed unit normal and  $S(\Sigma)$ the restriction of the spinor bundle $S(M)$ to $\Sigma$. Then 
\begin{align*}
  		 \rho_{\Sigma} : \Sol(D) & \longrightarrow\coinf(\Sigma, S(\Sigma)) \\
  		 \varphi & \longrightarrow \varphi_{\upharpoonleft \Sigma}
\end{align*}
is surjective.  Equivalently 
the Cauchy problem
\begin{equation*}
  \left\{
      \begin{aligned}
        D \psi & = 0 \\
        \rho_{\Sigma} \psi & = f, f \in\coinf(\Sigma; S(\Sigma))
      \end{aligned}
    \right.
\end{equation*}
is globally well-posed, the solution being denoted by $\psi = U_{\Sigma} f$.
We have:
\beq\label{cauchy}
U_{\Sigma} f = - \int_{\Sigma} G(x,y)\gamma(n(y))f(y)dVol_h,
\eeq
where $h$ is the Riemannian metric induced by $g$ on $\Sigma$.
We equip $\coinf(\Sigma,S(\Sigma))$ with the Hermitian form
\begin{equation} \label{eq3.1}
	(f_1|f_2)_{\Sigma} \defeq  \int_{\Sigma} \overline{f_1}\cdot \beta f_2 dVol_h.
\end{equation}

For $u\in \coinf(\Sigma;S(\Sigma))$, we define $\rho_{\Sigma}^*u \in \cD'(M; S(M))$ by
\[\int_M \overline{\rho_{\Sigma}^*u} \cdot \beta v dVol_g \defeq  \int_{\Sigma} \overline{u} \cdot \beta \rho_{\Sigma}v dVol_h, v \in \cinf(\Sigma;S(\Sigma)),\]
The identity \eqref{cauchy} can be rewritten as
\[ U_{\Sigma} f = (\rho_{\Sigma}G)^* \gamma(n)f, f\in\coinf(\Sigma;S(\Sigma)).\]

\subsubsection{Quantization of the Dirac equation}

For $\psi_1, \psi_2 \in \Sol(D)$ we set
\beq\label{defidef}
 \psi_1 \cdot \nu \psi_2 = \int_{\Sigma} \i J_{\mu}(\psi_1,\psi_2) n^{\mu} dVol_h = (\rho_{\Sigma} \psi_1 | \i\gamma(n) \rho_{\Sigma} \psi_2)_{\Sigma}.
 \eeq
Since $\nabla^{\mu}J_{\mu}(\psi_1,\psi_2)=0$, the right-hand side of \ref{eq3.1} is independent of the choice of $\Sigma$. From \eqref{e2.5} and the uniqueness of the Cauchy problem, we obtain that $\nu$ is positive definite.

Setting
\beq\label{nusig} \bar{f}_1 \cdot \nu_{\Sigma} f_2 \defeq  \int_{\Sigma} \i \bar{f}_1 \cdot \beta  \gamma(n) f_2 dVol_h, \ f_1,f_2 \in\coinf(\Sigma,S(\Sigma))
\eeq
we obtain that
\[\rho_{\Sigma} : (\Sol(D),\nu) \rightarrow (\coinf(\Sigma;S(\Sigma)),\nu_{\Sigma}), \]
is unitary, with inverse $U_{\Sigma}$.
We also get,  see e.g. \cite{Di}, that $G :\coinf(M; S(M)) \rightarrow \Sol(D)$ is surjective with kernel $D\coinf(M; S(M))$ and
\[ G : (\coinf(M; S(M))/ D\coinf(M; S(M)),\i(\cdot |G\cdot))  \rightarrow (\Sol,\nu) \]
is unitary.
Summarizing, the maps
\begin{equation} \label{unitarymaps}
	(\coinf(M; S(M))/ D\coinf(M; S(M)),\i(\cdot |G\cdot))  \xrightarrow{\, \, G \, \,} (\Sol,\nu)  \xrightarrow{\, \, \rho_{\Sigma} \, \,} (\cinf(\Sigma, S(\Sigma)),\nu_{\Sigma}),
\end{equation}
are unitary.

\subsection{Hadamard states for Dirac fields}\label{sec1b.3}

We denote by $\CAR(D)$ the $*$-algebra $\CAR(\cY,\nu)$ with $(\cY,\nu)$ being one of the equivalent pre-Hilbert spaces in (\ref{unitarymaps}).
We use the Hermitian form $(\cdot|\cdot)_{M}$ to pair $C_0^{\infty}(M,S(M))$ with $\cD'(M; S(M))$, and to identify continuous sesquilinear forms on $\coinf(M; S(M))$ with continuous linear maps from $\coinf(M; S(M))$ to $\cD'(M; S(M))$.
A quasi-free state $\omega$ on $\CAR(D)$ is defined by its {\em spacetime covariances} $\Lambda^{\pm}$, which satisfy
\beq\label{ploc}
\begin{array}{rl}
&\Lambda^{\pm} :\coinf(M,S(M)) \rightarrow \mathcal{D}'(M,S(M)) \text{ are linear continuous}, \\[2mm]
	&\Lambda^{\pm} \geq 0 \hbox{ for } (\cdot| \cdot)_{M}, \\[2mm]
	&\Lambda^+ + \Lambda^- = \i G,\\[2mm]
	&D \circ \Lambda^{\pm} = \Lambda^{\pm} \circ D = 0.
\end{array}	
\eeq
Alternatively, on can define $\omega$ by its {\em Cauchy surface covariances} $\lambda^{\pm}_{\Sigma}$, which satisfy
\beq\label{plic}
\begin{array}{rl}
	&\lambda^{\pm}_{\Sigma} :\coinf(\Sigma,S(\Sigma)) \rightarrow \mathcal{D}'(\Sigma,S(\Sigma))\text{ are linear continuous}, \\[2mm]
	&\lambda^{\pm}_{\Sigma} \geq 0  \hbox{ for } (\cdot| \cdot)_{\Sigma}, \\[2mm]
	&\lambda^+_{\Sigma} + \lambda^-_{\Sigma} = \i \gamma(n).
\end{array}	
\eeq
One can  go from one to the other, using the following identities which are easy to prove:
\beq\label{blurp}
\begin{array}{rl}
&\Lambda^{\pm} = (\rho_{\Sigma}G)^* \lambda^{\pm}_{\Sigma}(\rho_{\Sigma}G), \\[2mm]
&\lambda^{\pm}_{\Sigma} = (\rho_{\Sigma}^* \gamma(n))^* \Lambda^{\pm}(\rho_{\Sigma}^* \gamma(n)).
\end{array}	
\eeq

By the Schwartz kernel theorem, we can identify $\Lambda^{\pm}$ with distributional sections in $\cD'(M \times M; S(M) \boxtimes S(M))$, still denoted by $\Lambda^{\pm}$.
\subsubsection{Characteristic manifold}
The principal symbol $d(x,\xi)$ of $D$ is given by
\[
d(x,\xi) = \gamma(g^{-1}(x)\xi), \ (x, \xi)\in T^{*}M.
\]
The characteristic manifold  is  $\cN\defeq  \{(x, \xi)\in T^{*}M\setminus\zero: d(x, \xi)\hbox{ invertible}\}$. Since $d(x, \xi)^{2}= \xi\dual g^{-1}(x)\xi\one$, we have
\[
\cN= \{(x, \xi)\in T^{*}M\setminus\zero: \xi\dual g^{-1}(x)\xi=0\},
\]
so the characteristic manifold of $D$ is the same as the one of scalar Klein-Gordon operators.

The following notions are well-known and originate from \cite{DH}.

If $\Gamma\subset T^{*}M\times T^{*}M$ is a conic set, then one sets
\[
\Gamma'= \{((x_{1}, \xi_{1}), (x_{2}, -\xi_{2})): ((x_{1}, \xi_{1}), (x_{2}, \xi_{2}))\in \Gamma\}.
\]

The two connected components of $\cN$ are the {\em positive/negative energy shells}:
\[
\cN^{\pm}\defeq\{(x, \xi)\in \cN: \pm \xi\dual v>0\hbox{ for }v\in T_{x}M\hbox{ future directed}\}.
\]

We set
\[\cC= \{(X_1,X_2) \in \cN \times \cN : X_1 \sim X_{2}\}\subset T^{*}M\times T^{*}M,
\]
where we write $X_1 \sim X_2$ if $X_1,X_2 \in \cN$ and $X_1$ and $X_2$ lie on the same integral curve of the Hamiltonian vector field $H_p$  associated to $p(x, \xi)= \xi\dual g^{-1}(x)\xi$.

For $X_1 \sim X_2$, we write $X_1\succ X_2$, resp. $X_2\prec X_2$ if $X_1$ comes strictly after, resp. before $X_2$ with respect to the natural parameter on the integral curve of $H_p$ through $X_1$ and $X_2$.
We introduce the following subsets of $\cC$
\[
\begin{array}{l}	
\mathcal{C}_{\rm F} = \{(X_1,X_2) \in \cC: X_1\prec X_2 \} \\[2mm]
\mathcal{C}_{\overline {\rm F}} = \{(X_1,X_2) \in \cC: X_1\succ X_2 \}\\[2mm]
\mathcal{C}^{\pm}=  \{(X_1,X_2) \in \cC: X_1\in \cN^{\pm} \}. 
\end{array}
\]
\begin{definition}
$\omega$ is a Hadamard state if
\[ \WF(\Lambda^{\pm})' \subset\mathcal{C}^{\pm}.
\]
\end{definition}
The following lemma is certainly well-known.
\begin{lemma}\label{lemino}
 $\omega$ is a Hadamard state iff 
 \beq\label{blurg}
 \WF(\Lambda^{\pm})'\subset \cN^{\pm}\times \cN^{\pm}.
 \eeq
 \end{lemma}
 \proof Let $P= D^{2}$, which has a scalar principal symbol equal to $\xi\dual g^{-1}(x)\xi\one$, and let $E_{\rm ret/adv}$ its retarded/advanced inverses, $E= E_{\rm ret}- E_{\rm adv}$.  The arguments of \cite{DH} extend to show that
 \[
 \WF(E)'\subset \{(X_{1}, X_{2}): X_{1}\sim X_{2}\}.
 \]
  By uniqueness of $G_{\rm ret/adv}$ we have $G= D E$ hence (see e.g. \cite[Prop. A.7]{S}):
   \[
 \WF(G)'\subset \{(X_{1}, X_{2}): X_{1}\sim X_{2}\}.
 \]
  Since $\Lambda^{+}+ \Lambda^{-}= \i G$ this implies the lemma. \qed
\subsubsection{Feynman inverse}
Let $\omega$ a Hadamard state with spacetime covariances $\Lambda^{\pm}$.
Then
\[G_{\rm F} \defeq  \i^{-1} \Lambda^+ + G_{\rm adv} = -\i^{-1}\Lambda^-+G_{\rm ret}\]
is a Feynman inverse of $D$, meaning that
\[
DG_{\rm F} = G_{\rm F} D = \id, \WF(G_{\rm F})' = \Delta \cup \mathcal{C}_{\rm F},\]
where $\Delta= \{(X,X): X \in T^*M\setminus \zero\}$ is the diagonal.

\subsection{Hadamard condition for Cauchy surface covariances}
The following proposition gives a sufficient condition for the 
Cauchy surface covariances $\lambda^{\pm}_{\Sigma}$ to define a Hadamard state. 
\begin{proposition}\label{prop15.0a}
 Let \[
 \lambda_{\Sigma}^{\pm}\eqdef \i \gamma(n)c^{\pm}
 \]
  be the Cauchy surface covariances of a quasi-free state $\omega$. Assume that $c^{\pm}$ are  continuous from $\coinf(\Sigma;\cS_{\Sigma})$ to $\cinf(\Sigma;\cS_{\Sigma})$  and from $\cE'(\Sigma; \cS_{\Sigma})$ to $\cD'(\Sigma; \cS_{\Sigma})$, and that for some neighborhood $U$ of $\Sigma$ in $M$ we have
 \beq\label{troup.100}
\WF(U_{\Sigma}\circ c^{\pm})'\subset (\cN^{\pm}\cup\cF)\times T^{*}\Sigma, \hbox{ over }U\times \Sigma,
\eeq
where $\cF\subset T^{*}M$ is a conic set with $\cF\cap \cN= \emptyset$.  Then 
 $\omega$ is a Hadamard state.
 \end{proposition}
 
 \proof
 Let $\Lambda^{\pm}$ be the spacetime covariances of $\omega$.  By  \eqref{blurp} $\Lambda^{\pm}= \pm \i^{-1}U_{\Sigma}c^{\pm}\circ (\varrho_{\Sigma}G)$. Note that we are allowed to compose the kernels $U_{\Sigma}c^{\pm}$ and $\varrho_{\Sigma}G$ since $\varrho_{\Sigma}G: \coinf(M; S(M))\to \coinf(\Sigma; S(\Sigma))$.
 
We have seen above that  $\WF(G)'\subset \cN\times \cN$. Using also   \eqref{troup.100} and \cite[Thm. 8.2.14]{H} we obtain
  \beq\label{e9.00b}
 \WF(\Lambda^{\pm})'\subset (\cN^{\pm}\cup \cF)\times \cN \cup (\cN^{\pm}\cup \cF)\times \zero \hbox{ over }U\times M,
 \eeq
 where we recall that $\zero\subset T^{*}M$ is the zero section. 
 Since  $\Lambda^{\pm*}= \Lambda^{\pm}$  for the scalar product \eqref{burk} we obtain  that $(X, X')\in \WF(\Lambda^{\pm})'$ iff $(X', X)\in \WF(\Lambda^{\pm})'$. Using that $ \cF\cap \cN= \emptyset$, we then deduce from \eqref{e9.00b} that
\[
  \WF(\Lambda^{\pm})'\subset \cN^{\pm}\times \cN^{\pm} \ \cup \ \cN^{\pm}\times \zero\ \cup\  \zero\times \cN^{\pm}\hbox{ over }U\times M.
 \]
 Since $\Lambda^{+}- \Lambda^{-}= \i G$ and  using once more that $\WF(G)'\subset\cN\times \cN$ and the hypothesis $\cF\cap \cN= \emptyset$ this implies that 
 \[
 \WF(\Lambda^{\pm})'\cap (\cN^{\pm}\times \zero\ \cup \ \zero\times \cN^{\pm})= \emptyset,
 \]
 which proves   \eqref{blurg} over $U\times M$. To extend \eqref{blurg}  to $M\times M$ we use that  $D\Lambda^{\pm}=0$ and argue as in \cite[Lemma 6.5.5]{DH}.  \qed

\subsection{Conformal transformations}\label{sec1b.4}

From (\ref{econfo.2}) we get 
\[G = W\tilde{G}W^*,\]
and setting \[
 U : f\in\cinf_0(\Sigma,S(M)_{\Sigma}) \mapsto e^{\frac{1-n}{2}u}f,
 \] we obtain the following commutative diagram, with all arrows unitary :
\[
\begin{tikzcd}
&(\coinf(M; S(M))/ D\coinf(M; S(M)),\i(\cdot |G\cdot))  \arrow[r,"G"] \arrow[d,"W^*"] &(\Sol,\nu)  \arrow[r,"\rho_{\Sigma}"] \arrow[d,"W^{-1}"] &(\cinf(\Sigma, S(\Sigma)),\nu_{\Sigma}) \arrow[d,"U"] \\
&(\coinf(M; S(M))/ \tilde{D}\coinf(M; S(M)),\i(\cdot |\tilde{G}\cdot))  \arrow[r,"\tilde{G}"] &(\Sol,\tilde{\nu})  \arrow[r,"\rho_{\Sigma}"] &(\cinf(\Sigma, S(\Sigma)),\tilde{\nu_{\Sigma}})
\end{tikzcd}
\]

Let $\Lambda^{\pm}$ be the spacetime covariances of a quasi-free state $\omega$ for $D$.
Then 
\[\tilde{\Lambda}^{\pm} = e^{-\frac{n-1}{2}u} \Lambda^{\pm} e^{\frac{n+1}{2}u} \]
are the spacetime covariances of a quasi-free state $\tilde{\omega}$ for $\tilde{D}$, and
\[\tilde{\lambda}^{\pm}_{\Sigma} = (U^*)^{-1} \lambda^{\pm}_{\Sigma}U^{-1},\]
are the Cauchy surface covariances of a quasi-free state $\tilde{\omega}$, if $\lambda^{\pm}_{\Sigma}$ are the Cauchy surface covariances of a quasi-free state $\omega$.

Clearly the Hadamard condition is preserved by conformal transformations.
 \section{Dirac operators on spacetimes of bounded geometry}\label{sec2}

 In this section we study Dirac operators on spacetimes of bounded geometry.  We first recall some definitions about Riemannian manifolds of bounded geometry and their extensions to the Lorentzian case.  A useful result that we will prove is the fact that if a spacetime of bounded geometry admits a spin structure, then the spin and spinor bundles can be assumed to be of bounded geometry too.
 \subsection{Manifolds of bounded geometry}\label{sec2.1}
 \subsubsection{Riemannian manifolds of bounded geometry}\label{sec2.1.1}
  We recall that a Riemannian manifold $(M, \hat{g})$ is of {\em bounded geometry} if its injectivity radius $r_{\hat{g}}$ is strictly positive and $\nabla^{k}R$ are bounded tensors where $R$ is the Riemann curvature tensor and $\nabla$ the covariant derivative associated to $\hat{g}$, see  e.g. \cite[Appendix 1]{Sh}.
  
 Since we are interested in Lorentzian manifolds $(M, g)$   the role of a Riemannian metric $\hat{g}$ of bounded geometry will simply be to  define various function spaces, like spaces of bounded tensors, Sobolev spaces etc. 
  
  An equivalent definition (see e.g. \cite[Thm. 2.2]{GOW}) is as follows: let us denote by $B_{n}(0, 1)$ the unit ball in $\rr^{n}$, by $\delta$ the flat metric on $\rr^{n}$  and by $\BT^{p}_{q}(B_{n}(0, 1),\delta)$ the space of $(q,p)$ tensors on $B_{n}(0, 1)$ which are bounded on $B_{n}(0, 1)$ together with all their derivatives. 
  
  $\BT^{p}_{q}(B_{n}(0, 1),\delta)$ is a Fr\'echet space and one can hence define  a bounded family of $(q, p)$ tensors in $\BT^{p}_{q}(B_{n}(0, 1),\delta)$. For example a family $(f_{i})_{i\in I}$ of functions on $B_{n}(0, 1)$ is bounded if $\sup_{i\in I, x\in B_{n}(0, 1)}| \p_{x}^{\alpha}f_{i}(x)|<\infty$ for all $\alpha\in \nn^{n}$.

\begin{proposition}\label{prop1.1}
 $(M, \hat{g})$ is of bounded geometry if and only if for each $x\in M$ there exists $U_{x}\subset M$ open neighborhood of $x$ and
\[
\psi_{x}: U_{x} \xrightarrow{\sim} B_{n}(0,1)
\]
a  diffeomorphism with $\psi_{x}(x)=0$ such that if $\hat{g}_{x}\defeq   (\psi_{x}^{-1})^{*}\hat{g}$ then:

\ben
\item  the family $\{\hat{g}_{x}\}_{x\in M}$ is  bounded in ${\rm BT}^{0}_{2}(B_{n}(0,1), \delta)$,

\item there exists $c>0$ such that:
\[
c^{-1}\delta\leq \hat{g}_{x}\leq c \delta, \ x\in M.
\]
\een
\end{proposition}
It is known that one can find  a sequence $(x_{i})_{i\in \nn}$ of points in $M$ such that if $(U_{i}, \psi_{i})\defeq (U_{x_{i}}, \psi_{x_{i}})$ then $(U_{i}, \psi_{i})_{i\in \nn}$ is an atlas of $M$ with the additional property that there exists $N\in \nn$ such that $\bigcap_{i\in J}U_{i}= \emptyset$ if $\sharp J>N$. Such atlases are called {\em bounded atlases} of $M$. 

The standard choice for  $(U_{i}, \psi_{i})$ is $(B^{\hat{g}}(x_{i}, r), \exp^{\hat{g}}_{x_{i}})$ for $0<r<r_{\hat{g}}$ i.e. the geodesic ball and exponential map at $x_{i}$. It follows that  without loss of generality we can assume that the $U_{i}$ are simply connected.

One can associate to a bounded atlas  $(U_{i}, \psi_{i})_{i\in \nn}$ a partition of unity
\[
1= \sum_{i\in \nn}\chi_{i}^{2}, \ \chi_{i}\in \coinf(U_{i}),
\]
such that $\{(\psi_{i}^{-1})^{*}\chi_{i}\}_{i\in \nn}$ is a bounded sequence in $\cinf_{\rm b}(B_{n}(0,1))$. Such a partition of unity is  called a {\em bounded partition of unity}.
\subsubsection{Bounded tensors, bounded differential operators, Sobolev spaces}\label{sec2.1.2}
We recall now several definitions due to Shubin \cite{Sh}, see also \cite[Subsect. 2.3]{GOW}.

If $(M, \hat{g})$ is a Riemannian manifold of bounded geometry,  we denote by  $\BT^{p}_{q}(M, \hat{g})$ the space of {\em bounded} $(q, p)$ tensors on $M$. 

If $I\subset \rr$ is an interval  and $\cF$ is a Fréchet space whose topology is defined by the seminorms $\| \cdot \|_{n}$ $n\in \nn$, we denote by $\cinfb(I, \cF)$ the space of maps $f: I\to \cF$ such that $\sup_{t\in I}\| \p_{t}^{p}f(t)\|_{n}<\infty$ for all $n, p\in \nn$. Equipped with the obvious seminorms it it itself a Fréchet space.

 We will use this convention to define the spaces $\cinfb(I, \BT^{p}_{q}(M, \hat{g}))$.

The space  ${\rm Diff}(M, \hat{g})$ is the space of bounded differential operators on $M$, i.e. differential operators which form a bounded family of differential operators on $B_{n}(0, 1)$ when expressed in a bounded atlas of $M$. 

Finally  we denote by  $H^{s}(M, \hat{g})$ the Sobolev space of order $s\in \rr$.

\subsubsection{Vector bundles of bounded geometry}\label{sec2.1.3}
Let $(M, g)$ be a Riemannian manifold of bounded geometry.    We recall  the definition of vector bundles of bounded geometry, see \cite{Sh}.  

A vector bundle $E\xrightarrow{\pi}M$ of rank $N$  is of {\em bounded geometry} if there exists a bounded covering $(U_{i})_{i\in \nn}$ of $M$  which form a bundle atlas of $E$ such that the transition maps $t_{ij}: U_{ij}\to M_{N}(\cc)$ are a bounded family of matrices. 

The space of bounded sections of $E$ (analog to the spaces of bounded tensors) is denoted by $\cinfb(M; E)$.  Concretely if $u$ is a section of $E$, we denote by $u_{i}: U_{i}\to \cc^{N}$ its local trivializations over $U_{i}$ and $u\in \cinfb(M; E)$ iff the family $(u_{i})_{i\in \nn}$ is bounded in $\cinfb(U_{i}; \cc^{N})$.

The notion of bounded differential operators acting on smooth sections of  $E$ is now defined as before, as is the notion of bounded Hermitian forms on (the fibers of) $E$.  
\subsubsection{Bounded Hilbert space structures}
Using a partition of unity one can equip a vector bundle $E\xrightarrow{\pi}M$  of bounded geometry with a positive definite bounded Hermitian form $\beta$, (ie a Hilbert space structure on the fibers of $E$).  This means that  if $\beta_{i}: U_{i}\to L_{\rm h}(\cc^{N}, \cc^{N*})$ are its local trivializations, then $\beta_{i}>0$ and the families $(\beta_{i})_{i\in \nn}$ and $(\beta_{i}^{-1})_{i\in \nn}$ are bounded in $\cinfb(U_{i}; L_{\rm h}(\cc^{N}, \cc^{N*}))$ and $\cinfb(U_{i}; L_{\rm h}(\cc^{N*}, \cc^{N}))$.

Any two of these bounded Hilbert space structures are equivalent, in the sense of the bounded geometry.

\subsubsection{Sobolev spaces}
One defines the Sobolev spaces $H^{s}(M; E)$ in the natural way, for example using the norm
\[
\| u\|_{s}^{2}= \sum_{i\in \nn}\| T_{i}\circ \psi_{i}\chi_{i}u\|^{2}_{s},
\]
where $(U_{i}, \psi_{i})_{i\in \nn}$ is a bounded atlas, $T_{i}: \pi^{-1}(U_{i}) \to U_{i}\times \cc^{N}$ are local trivialisations, $1= \sum \chi_{i}^{2}$ is a bounded partition of unity subordinate to $(U_{i})_{i\in \nn}$ and $\| \cdot\|_{s}$ is the usual Sobolev norm on $\rr^{n}$.

The topology of $H^{s}(M; E)$ is independent of the above choices.

\subsubsection{Principal bundles of bounded geometry}\label{sec2.1.4}
Similarly one can define principal bundles of bounded geometry.
\begin{definition}\label{def1.10}
 Let $G\subset M_{N}(\cc)$ be  a matrix Lie group. A $G$-principal bundle $G\to P\xrightarrow{\pi}M$  is  {\em  of bounded geometry}  if there exists a bounded covering $(U_{i})_{i\in \nn}$ of $M$  which form a bundle atlas of $P$ such that the transition maps $t_{ij}: U_{ij}\to M_{N}(\cc)$ are a bounded family of matrices. 
\end{definition}
Clearly if $G\to P\xrightarrow{\pi}M$ is a principal bundle of bounded geometry   and  $\rho: G\to M_{N}(\cc)$ is a matrix representation of $G$ then the associated bundle $P\times_{\rho} \cc^{N}$ is a vector bundle of bounded geometry.

\subsubsection{Bounded frames}\label{sec2.1.5}
\begin{definition}\label{def2.1}
  If $(U_{i},\psi_{i})_{i\in \nn}$ is a bounded  atlas, a family of local frames $(\mathcal{F}_{i})_{i\in \nn}$ of $TM$ over $U_{i}$ is called a bounded family of frames if $\mathcal{F}_{i}= (f_{i, a})_{0\leq a\leq d}$ and the families $(f_{i, a})_{i\in \nn}$ are bounded.
\end{definition}
 Of course a manifold of bounded geometry admits bounded families of frames, for example the frames associated to local normal coordinates for the reference metric $\hat{g}$.
 
 \subsection{Lorentzian manifolds of bounded geometry}\label{sec2.2}
 We now recall the notion of a Lorentzian manifold $(M, g)$ of bounded geometry, with respect to a reference Riemannian metric $\hat{g}$, see \cite[Sect. 3]{GOW}.
 \begin{definition}\label{defp3.1}
 Let $M$ a manifold equipped with a reference Riemannian metric $\hat{g}$ of bounded geometry. A Lorentzian metric $g$ on $M$ is of {\em bounded geometry} if $g\in {\rm BT}^{0}_{2}(M, \hat{g})$ and $g^{-1}\in {\rm BT}^{2}_{0}(M, \hat{g})$.
\end{definition}
 \subsubsection{Bounded orthonormal frames}\label{sec2.2.1}
\begin{proposition}\label{prop1.2}
 Let $(M, g)$ a spacetime of bounded geometry, with respect to the reference Riemannian metric $\hat{g}$.  Then  $(M, g)$ admits a bounded family of local orthonormal oriented and time oriented frames.
\end{proposition}
\proof
  Let  $(\mathcal{F}_{i})_{i\in \nn}$  a bounded family of frames, $\mathcal{F}_{i}= (f_{i, a})_{0\leq a\leq d}$, and 
$\pmb{g}_{i, ab}= f_{i, a}\dual g f_{i, b}$.  Without loss of generality we can assume that $f_{i, 0}$ is future directed. The family of matrices $(\pmb{g}_{i})_{i\in \nn}$ is bounded. 

We choose a neighborhood $V_{0}$ of $\eta$ as in Lemma  \ref{lemma-app.1} small enough   such that  $\pmb{t}f_{i, 0}$ is future directed for all $\pmb{t}\in U_{0}= F(V_{0})$. We also choose 
  $U_{i}= B^{g}(x_{i}, r)$ for $0<r\ll 1$  such that $\pmb{g}_{i}\in V_{0}$.

 We set then $\pmb{t}_{i}\defeq  F(\pmb{g}_{i})$, which is a bounded family of matrices, $\mathcal{E}_{i}\defeq  \mathcal{F}_{i}\pmb{t}_{i}$ and $(\mathcal{E}_{i})_{i\in \nn}$ is then a bounded family of orthonormal oriented and time oriented frames.\qed
 \subsection{Cauchy surfaces of bounded geometry}\label{sec2.3}
 Let $(M, g)$ a spacetime of bounded geometry, with respect to a reference Riemannian metric $\hat{g}$. If $(M, g)$ is globally hyperbolic, then $M$ admits Cauchy surfaces. We recall the notion of {\em Cauchy surfaces of bounded geometry} introduced in \cite{GOW}.
\begin{definition}\label{defp3.3}
 Let $(M, g)$ be a  Lorentzian manifold of bounded geometry with respect to the reference Riemannian metric  $\hat{g}$. Assume  that  $(M, g)$ is  globally hyperbolic and let $\Sigma\subset M$ a smooth {\em spacelike} Cauchy hypersurface.   Then $\Sigma$ is called a {\em bounded geometry Cauchy hypersurface} if:
 \ben
 \item  the embedding  $i: \Sigma\to M$ is of bounded geometry for $\hat{g}$,
 \item if $n(y)$ for  $y\in \Sigma$ is the future directed unit normal  for $g$ to $\Sigma$ one has:
 \[
\sup_{y\in \Sigma}n(y)\cdot \hat{g}(y)n(y)<\infty.
\]
 \een
\end{definition} We refer to \cite[Def. 2.8]{GOW} for the notion of an embedding of bounded geometry, due to \cite{eldering}.

If $\Sigma$ is a spacelike Cauchy surface of bounded geometry, then the Gaussian normal coordinates to $\Sigma$ are well adapted to the bounded geometry framework. We recall a  result in this direction, see \cite[Thm. 3.5]{GOW}.
\begin{theoreme}\label{th-omar}
 Let $(M, g)$ a Lorentzian manifold of bounded geometry and $\Sigma$ a bounded geometry Cauchy hypersurface. Then the following holds:
\ben
\item there exists $\delta>0$ such that the normal geodesic flow to $\Sigma$:
\[
\chi:\ \begin{array}{l}
]-\delta, \delta[\times \Sigma\to M\\[2mm]
(s, y)\mapsto \exp_{y}^{g}(sn(y))
\end{array}
\] 
is well defined and is a smooth diffeomorphism on its image;
\item $\chi^{*}g= - ds^{2}+h_{s}$, where $\{h_{s}\}_{s\in\,]-\delta, \delta[}$ is a smooth family of Riemannian metrics on $\Sigma$ with
\[
\begin{array}{rl}
i)&(\Sigma, h_{0})\hbox{ is of bounded geometry},\\[2mm]
ii)& s\mapsto h_{s}\in \cinf_{\rm b}(]-\delta, \delta[, \BT^{0}_{2}(\Sigma, h_{0})), \\[2mm]
iii)& s\mapsto h^{-1}_{s}\in \cinf_{\rm b}(]-\delta, \delta[, \BT^{2}_{0}(\Sigma, h_{0})).
\end{array}
\]
 \een
\end{theoreme}
We recall that the spaces $\cinfb(I, \BT^{p}_{q}(M, \hat{g}))$ are defined in \ref{sec2.1.2}.
\subsection{Spin structures on manifolds of bounded geometry}\label{sec2.4}
Let $(M, g)$ a spacetime of bounded geometry with respect to a reference Riemannian metric $\hat{g}$.  By Prop. \ref{prop1.2} we know that the principal bundle $P\SO^{\uparrow}(M, g)$ is of bounded geometry. Therefore the associated Clifford bundle $\Cl(M, g)$ is of bounded geometry.

Assume that $(M, g)$ admits a spin structure and let $P\Spin(M , g)$ the spin bundle over $M$.  By Prop. \ref{propadd1.1} if $\pmb{o}_{ij}: U_{ij}\to \SO^{\uparrow}(1, d)$ are the transition maps of $P\SO^{\uparrow}(M, g)$, we can assume that the transition maps $\pmb{s}_{ij}$ of $P\Spin(M , g)$ satisfy $Ad(\pmb{s}_{ij})= \pmb{o}_{ij}$.  

Since $Ad: \Spin^{\uparrow}(1, d)\to \SO^{\uparrow}(1, d)$ is a two sheeted covering, the
transition maps $\pmb{s}_{ij}$ form a bounded family of matrices and hence $P\Spin(M , g)$  is a principal bundle of bounded geometry. 

The associated bundle $S(M)$ is then a vector bundle of bounded geometry. If the  section $m\in \cinf(M; L(S(M)))$ is a bounded section, then the Dirac operator $D$ defined in Def. \ref{def1.2} is a bounded differential operator.

\subsubsection{The case of cartesian products}\label{carproduct}
We  now precise the results of Subsect. \ref{spinprod} in the bounded geometry case.

Let  hence $I= ]-\delta, \delta[$, $(\Sigma, k)$ be  a Riemannian manifold of bounded geometry  and $M = I\times \Sigma$. Then setting $\hat{g}= dt^{2}+ k$, $(M, \hat{g})$ is of bounded geometry. If $(V_{i})_{i\in \nn}$ is a bounded covering for $(\Sigma, k)$ and $U_{i}= I\times V_{i}$ then $(U_{i})_{i\in \nn}$ is a bounded covering for $(M, \hat{g})$.

Let \[
g= - dt^{2}+ h_{t}(\rx)d\rx^{2}
\]
 be such that $g$ is of bounded geometry with respect to $\hat{g}$, or equivalently that $I\ni t\mapsto h_{t}$ satisfies the properties in (2) of Thm. \ref{th-omar}, (with $h_{0}$ replaced by $k$).  
 
 Since $h_{0}$ and $k$ are equivalent as Riemannian metrics of bounded geometry (see \cite[Subsect. 2.5]{GOW} for the terminology), we can assume that $k= h_{0}$. 
 
 The induced spin structure on $(\Sigma, h_{0})$ is of bounded geometry and the restriction $S_{t}(\Sigma)= S(\Sigma)$ of $S(M)$ to $\Sigma_{t}= \{t\}\times \Sigma$ is a vector bundle of bounded geometry.

   \subsection{Hypotheses}\label{sec5.0}
   In this subsection we describe the geometric framework introduced in    \cite[Sect. 3.3]{GOW} that we will use in Sect. \ref{sec5} to construct Hadamard states.  Roughly speaking we consider spacetimes of bounded geometry with a bounded geometry Cauchy surface. 
   
   Because of the covariance of the Dirac equation under conformal transformations, our construction applies actually to a wider class of spacetimes. Let us now describe more precisely our framework.

If $(M_{i}, g_{i})$, $i=1, 2$ are two space\-times, a \emph{space\-time embedding} $i: (M_{1}, g_{1})\to (M_{2}, g_{2})$ is  an isometric embedding preserving  the orientation and time-orientation. In addition, if $(M_{i}, g_{i})$ are globally hyperbolic, one says that $i$  is {\em causally compatible} if:
\[
I^{\pm, g_{1}}_{M_{1}}(U)= i^{-1}(I^{\pm, g_{2}}_{M_{2}}(i(U)), \ \forall \, U\subset M_{1},
\]
where $I^{\pm,g}_{M}(K)$ denote the future/past causal shadows of $K\subset M$ for the metric $g$.

We fix an even dimensional  globally hyperbolic space\-time $(M,g)$ admitting a spin structure, a smooth space-like Cauchy hypersurface $\Sigma$,  and a function $m\in \cinf(M; L(S(M)))$.
We assume that there exist:
\ben
\item  a neighborhood $U$ of $\Sigma$ in $M$, 
\item a Lorentzian metric $\tilde{g}$ on $M$,
\item a function $\tilde{u}\in \cinf(M; \rr)$,
\een
such that:

(H1)   $(M, e^{2\tilde{u}}\tilde{g})$ is globally hyperbolic, $Id: (U, g)\to (M, e^{2\tilde{u}}\tilde{g})$ is causally compatible,

(H2)  $\tilde{g}$ is of bounded geometry for some reference Riemannian metric $\hat{g}$, $\Sigma$ is a Cauchy hypersurface of  bounded geometry in $(M, \tilde{g})$,

(H3) $d\tilde{u}$ belongs to $\BT^{0}_{1}(M, \hat{g})$,

(M) $e^{-\tilde{u}}m$ belongs to $\BT^{0}_{0}(M, L(S(M)))$.

The following result is proved in \cite[Prop. 3.7]{GOW}.
\begin{proposition}
	Assume hypotheses {\rm (H)}. Then there exist:
	\begin{enumerate}
	\item an open interval $I$ with $0 \in I$, a diffeomorphism $\chi : I \times \Sigma \to U$,
	\item a smooth family $\{h_t\}_{t\in I}$ of Riemannian metrics on $\Sigma$ with
	
		$(\Sigma,h_0)$ is of bounded geometry,
		
		$t \in I \to h_t \in \cinfb(I;BT^0_2(\Sigma,h_0))$, $t \in I \to h_t^{-1} \in \cinfb(I;BT^2_0(\Sigma,h_0))$
	\item a function $u \in \cinf(I\times \Sigma)$,
		
		$du \in BT^0_1(I\times \Sigma,dt^2+h_0)$,
	\end{enumerate}
such that
\[ \chi^*g = e^{2u(t,y)}(-dt^2+h_t(y)dy^2) \text{ on $I \times \Sigma$.} \]
If moreover hypothesis {\rm (M)} holds then
\[e^{-u}m \circ \chi^{-1} \in \cinfb(I;BT^0_0(\Sigma,h_0)).\]
\end{proposition}
\subsubsection{Examples}
Many  spacetimes of physical interest satisfy hypotheses {\rm (H)}, see  \cite[Sect. 4]{GOW} for more details. For example cosmological spacetimes, Kerr and Kerr-de Sitter  exterior regions,  the  Kerr-Kruskal extension,  and double cones, wedges, and future past lightcones in Minkowski spacetime satisfy {\rm (H)}. In all these examples a constant mass term $m$ satisfy hypothesis {\rm (M)}.

All the examples above  are parallelizable, hence admit a unique spin structure.

 \section{Pseudodifferential calculus on manifolds of bounded geometry}\init\label{sec4}

 We refer the reader to \cite{Sh}, \cite{Ko}.
 We fix a manifold $M$ equipped with a reference Riemannian metric $\hat{g}$ of bounded geometry and a complex vector bundle  $E\xrightarrow{\pi}M$ of bounded geometry.   Shubin's calculus, introduced in \cite{Sh}, is a global pseudodifferential calculus on $M$, adapted to the bounded geometry defined by $\hat{g}$. It generalizes both the pseudodifferential calculus on a compact manifold and the uniform pseudodifferential calculus on $\rr^{n}$. 
 
 An important result is  {\em Seeley's theorem} proved in \cite{alnv}, which states that if $A$ is an elliptic, positive selfadjoint pseudodifferential operator of order $m$ on $M$, then its complex powers $A^{z}$ form a holomorphic family of pseudodifferential operators of order $m\Re z$.
 
 The framework generalizes easily to 'matrix valued' pseudodifferential operators, acting on sections of a bounded vector bundle. Many results go through without changes, but some of them need more care and will be treated in some details in this section.
 
 We will also consider {\em time dependent} pseudodifferential operators $A= A(t)$ depending on $t\in I$ an open interval. The framework of \cite{alnv} is general enough to accomodate this extension without much additional work.

 \subsection{Symbol classes}\label{sec4.1}
 \subsubsection{Symbol classes on $\rr^{n}$}
 If $U\subset \rr^{n}$ is an open set, 
 we denote by $S^{m}(T^{*}U; L(\cc^{N}))$ the set of functions $a\in\cinf(T^{*}U; L(\cc^{N}))$ such that
 \[
\sup_{(x, \xi)\in T^{*}U} |\langle \xi\rangle^{- m+ |\beta|}\p_{x}^{\alpha}\p_{\xi}^{\beta}a(x, \xi)|<\infty, \ \forall \alpha, \beta\in \nn^{n},
 \]
 equipped with its Fr\'echet space topology. 

We  denote by $E$  a continuous extension
\[
E: S^{m}(T^{*}B(0, 1); L(\cc^{N}))\to S^{m}(T^{*}\rr^{n}; L(\cc^{N})).
\]

 \subsubsection{Symbol classes on $M$}\label{sec4.1.1}
 We fix $(U_{i}, \psi_{i})_{i\in \nn}$ a bounded atlas of $M$, $t_{i}: \pi^{-1}(U_{i})\tosim U_{i}\times \cc^{N}$ local trivialisations of $E$ and $1= \sum \chi_{i}^{2}$ a bounded partition of unity subordinate to $\{U_{i}\}_{i\in \nn}$. 
 
If $u\in \cinf(U_{i}; E)$ resp. $a\in \cinf(T^{*}U_{i}; L(E))$ we  denote by $T_{i}u\in \cinf(B(0, 1); \cc^{N})$ resp. $\tilde{T}_{i}a\in \cinf(T^{*}B_{n}(0, 1); L(\cc^{N}))$ the pushforward of $u$, resp. $a$
 obtained from $\psi_{i}$ and $t_{i}$.  
\begin{definition}\label{def4.1}
 For $m\in \rr$ we denote by $S^{m}(T^{*}M; L(E))$ the set of sections $a\in \cinf(T^{*}M; End(E))$ such that $\{\tilde{T}_{i}a\}_{i\in \nn}$ is a bounded family in 
  $S^{m}(T^{*}B(0, 1); L(\cc^{N}))$.
\end{definition}
 We denote by   $S^{m}_{\rm ph}(T^{*}M; L(E))$ the subspace of poly-homogeneous symbols, which is defined as in the scalar case, see e.g. \cite[Subsect. 5.2]{GOW}.  The space of homogeneous symbols of order $m$ is denoted by $S^{m}_{\rm h}(T^{*}M; L(E))$. 
 
 The spaces $S^{m}_{\rm (ph)}(T^{*}M; L(E))$ are equipped with Fr\'echet space topologies (although the topology of $S^{m}_{\rm ph}(T^{*}M; L(E))$
 deserves some care, see \cite{alnv} or \cite[ 5.2.1]{GOW}).
 
To ease notation these spaces will often be simply denoted by $S^{m}(T^{*}M), S^{m}_{\rm ph}(T^{*}M)$, or even $S^{m}$, $S^{m}_{\rm ph}$.
 
 \subsubsection{Time dependent symbols}\label{sec4.1.2}
The notation $\cinfb(I; \cF)$ for $I$ an interval and $\cF$ a Fréchet space was defined in \ref{sec2.1.2}.
We use this notation to define the spaces $\cinfb(I; S^{m}_{({\rm ph})}(T^{*}M; L(E)))$. For example $\cinfb(I; S^{m}(T^{*}\rr^{n}))$ is the space of smooth functions $a: I\times T^{*}\rr^{n}\to \cc$ such that 
\[
\sup_{I\times T^{*}\rr^{n}}\langle \xi\rangle^{-m+|\beta|}|\p_{t}^{k}\p_{x}^{\alpha}\p_{\xi}^{\beta}a(t, x, \xi)|<\infty, \ k\in\nn, \alpha, \beta\in \nn^{n}.
\]

\subsection{Pseudodifferential operators}\label{sec4.2}

 Let us fix a bounded Hilbertian structure $(\cdot| \cdot)_{E}$ on the fibers of $E$ and  define the scalar product
 \[
 (u|v)= \int_{M}(u(x)| v(x))_{E}dVol_{g}, \ u, b\in \coinf(M; E).
 \]
We denote by $L^{2}(M; E)$ the completion of $\coinf(M; E)$ for the induced norm.  Different bounded Hilbertian structures produce equivalent norms. 
\subsubsection{Quantization map}\label{sec4.2.1}
If $a\in S^{m}_{\rm ph}(T^{*}M; L(E)))$ we set
\[
\Op(a)u\defeq  \sum_{i\in \nn}\chi_{i} T_{i}^{-1}\circ \Op (E \tilde{T}_{i}a)\circ T_{i} \chi_{i}u, \ u\in \coinf(M; E),
\]
where the maps $T_{i}, \tilde{T}_{i}$ were defined in Subsect. \ref{sec4.1}, $1= \sum_{i}\chi_{i}^{2}$ is a bounded partition of unity subordinate to $(U_{i})_{i\in \nn}$ and the map $\Op$ is the usual Kohn-Nirenberg quantization:
\[
\Op(a)u(x)= (2\pi)^{-n}\int\e^{\i (x-y)\cdot \xi}a(x, \xi)u(y)dyd\xi, 
\]
 for $a\in S^{m}(T^{*}\rr^{n}; L(\cc^{N}))$, $u\in\coinf(\rr^{n}; \cc^{N})$.

We use the same definition for $a(t)\in \cinfb(I; S^{m}_{({\rm ph})}(T^{*}M; L(E)))$.

\subsubsection{Ideals of smoothing operators}\label{sec4.2.1b}

The quantization map $\Op$ depends of course of the choice of the atlas $(U_{i}, \psi_{i})_{i\in \nn}$, local trivializations $t_{i}$, cutoff functions $\chi_{i}$ and the extension map $E$.  If we denote by $\Op'$ another quantization map for a different choice of the above data, one has
\[
\Op(a)- \Op'(a)\in \cW^{-\infty}(M; L(E)),
\]
where
 \[
\cW^{-\infty}(M; L(E))\defeq   \bigcap_{m\in \nn}B(H^{-m}(M; E), H^{m}(M; E)),
\]
equipped with the topology given by the seminorms 
\[
\|A\|_{m}=\|(- \Delta_{\hat{g}}+1)^{m/2} A(- \Delta_{\hat{g}}+1)^{m/2}\|_{B(L^{2}(M; E))},
\]
$-\Delta_{\hat{g}}$ being the (scalar) Laplace-Beltrami operator on $(M; \hat{g})$. We will often denote $\cW^{-\infty}(M; L(E))$ simply by $\cW^{-\infty}(M)$ or even $\cW^{-\infty}.$

For time dependent symbols we have similarly
\[
\Op(a(t))- \Op'(a(t))\in\cinfb(I; \cW^{-\infty}(M; L(E))),
\]
where the later space is equipped with the
 topology given by the seminorms
\[
\|A\|_{m,p}= \sup_{t\in I, k\leq p}\|\p_{t}^{k}A(t)\|_{m}.
\]
An important property of $\cW^{-\infty}$ and $\cinfb(I; \cW^{-\infty})$ is the so-called {\em spectral invariance} property  introduced in \cite{alnv}.
It allows to 'close' the pseudodifferential calculus with respect to the operation of taking inverses. For completeness we state this property. Its proof is exactly the same as in the scalar case, see \cite[Lemma 5.5]{GOW}.
\begin{lemma}\label{spectral1}
 Let $R_{-\infty}\in \cW^{-\infty}(M; L(E))$ such that $\one - R_{-\infty}$ is invertible in $B(L^{2}(M; E))$. Then 
 \[
(\one - R_{-\infty})^{-1}= \one -R_{1, -\infty}\hbox{ for }R_{1, -\infty}\in \cW^{-\infty}(M; L(E)).
\]
\end{lemma}
The same result holds replacing $\cW^{-\infty}(M; L(E))$ by $\cinfb(I; \cW^{-\infty}(M; L(E)))$ and $L^{2}(M; E)$ by $L^{2}(I; L^{2}(M; E))$.
\begin{lemma}\label{spectral2}
 Let $R_{-\infty}\in \cinfb(I;\cW^{-\infty}(M; L(E)))$ such that $\one - R_{-\infty}$ is invertible in $B(L^{2}(I;L^{2}(M; E))$. Then 
 \[
(\one - R_{-\infty})^{-1}= \one -R_{1, -\infty}\hbox{ for }R_{1, -\infty}\in \cinfb(I;\cW^{-\infty}(M; L(E))).
\]
\end{lemma}

\subsubsection{Pseudodifferential operators}\label{sec4.2.2}

We equip $\Op(S_{\rm ph}^{m}(T^{*}M; L(E)))$ by the Fr\'echet space topology inherited from $S^{m}_{\rm ph}(M; L(E))$. 

Note that $\Op(S_{\rm ph}^{m}(T^{*}M; L(E)))$ is not stable under composition. To obtain an algebra of operators, one has to add the ideal $\cW^{-\infty}(M; L(E))$.
\begin{definition}\label{def4.2}
  For $m\in \rr$ we set 
  \[
  \Psi^{m}(M; L(E))\defeq  \Op(S^{m}_{\rm ph}(T^{*}M; L(E)))+ \cW^{-\infty}(M; L(E)),
  \]
We equip $ \Psi^{m}(M; L(E))$ with the Fr\'echet space topology obtained from those of $\Op(S_{\rm ph}^{m}(T^{*}M; L(E)))$ and $\cW^{-\infty}(M; L(E))$, the decomposition of $A\in \Psi^{m}$ as $A= \Op(a)+ R_{-\infty}$ being unique, once the quantization map is fixed. 
\end{definition}
We define similarly the spaces of time dependent pseudodifferential operators:
\[
\cinfb(I;  \Psi^{m}(M; L(E)))\defeq  \cinfb(I; \Op(S^{m}_{\rm ph}(T^{*}M; L(E))))+ \cinfb(I; \cW^{-\infty}(M; L(E))).
\]
\subsubsection{Principal symbol}\label{sec4.2.3}
\begin{definition}\label{def4.3}
 Let $A= \Op(a)+R_{-\infty}\in \Psi^{m}(M; L(E)$.  We denote by $\sigma_{\rm pr}(A)$ the {\em principal symbol} of $A$ defined as 
 \[
 \sigma_{\rm pr}(A)\defeq  [a]\in S^{m}_{\rm ph}(T^{*}M; L(E))/S^{m-1}_{\rm ph}(T^{*}M; L(E)).
 \] 
  $\sigma_{\rm pr}(A)$ is independent on the decomposition of $A$ as $\Op(a)+R_{-\infty}$ and on the choice of the good quantization map $\Op$.
\end{definition}
We use the analogous definition for $A\in \cinfb(I; \Psi^{m}(M, L(E)))$, the principal symbol $\sigma_{\rm pr}(A)$ being now in $\cinfb(I; S^{m}_{\rm h}(T^{*}M; L(E)))$. 

As usual we choose a representative of  $ \sigma_{\rm pr}(A)$ which is homogeneous of order $m$ on the fibers of $T^{*}M$ so we can assume that $\sigma_{\rm pr}(A)\in S^{m}_{\rm h}(T^{*}M; L(E))$.
\subsubsection{Ellipticity}\label{sec4.2.4}
\begin{definition}
An operator $A\in\Psi^{m}(M; L(E))$   is {\em elliptic} if  $\sigma_{\rm pr}(A)(x, \xi)$  is invertible for all $(x, \xi)\in T^{*}M$  and
 \beq\label{e4.1}
\sup_{ (x, \xi)\in T^{*}M, | \xi|= 1}\|\sigma_{\rm pr}(A)^{-1}(x, \xi)\|<\infty.
\eeq
\end{definition}
To define the norm in \eqref{e4.1} we choose a bounded Hilbertian structure on the fibers of $E$, the definition being independent on its choice.

We use the same definition for time dependent operators $A\in\cinfb(I; \Psi^{m}(M; L(E)))$, requiring also uniformity w.r.t. $t\in I$ in 
\eqref{e4.1}.

An elliptic operator $A\in\Psi^{m}(M; L(E))$  has a parametrix $B\in \Psi^{-m}(M; L(E))$, unique modulo $\cW^{-\infty}(M; L(E))$ such that $AB- \one, BA- \one\in\cW^{-\infty}(M; L(E))$. A parametrix of $A$ will be denoted by $A^{(-1)}$.
 
 The following result is well-known in  the scalar case, using the existence of a parametrix, see \cite{Ko}.
 \begin{proposition}\label{prop4.1}
Let $A\in \Psi^{m}(M; L(E))$, $m\geq 0$ elliptic. Then the following holds:
\ben
\item  $A$ is closeable on $\coinf(M; E)$ with $\Dom A^{\rm cl}= H^{m}(M; L(E))$.  
\item  if  $0\not\in \sigma(A^{\rm cl})$, then $A^{-1}\in \Psi^{-m}(M; L(E))$ and
\[
\sigma_{\rm pr}(A^{-1})= (\sigma_{\rm pr}(A))^{-1}.
\]
\een
\end{proposition}

\proof (1) follows easily from the existence of a parametrix $A^{(-1)}$. To prove (2) we first reduce ourselves to the case $m=0$. In fact if $Q= (-\Delta_{\hat{g}}+ 1)^{-m/2}$, we know that $Q\in \Psi^{-m}$ since Prop. \ref{prop4.1} is true in the scalar case (see for example \cite{Ko}), and $QP\in \Psi^{0}$ is invertible. 

If $P\in \Psi^{0}$ is invertible,  we choose a sequence $f_{n}$ of polynomials such that $f_{n}(\lambda)\to \lambda^{-1}$, uniformly on $[a, b]$ where  $\sigma(P^{*}P)\subset [a, b]$.  Then $Q_{n}= f_{n}(P^{*}P)P^{*}\in \Psi^{0}$ and $Q_{n}\to P^{-1}$ in $B(L^{2}(M; E))$. 
 
 Let $Q\in \Psi^{0}$ a parametrix of $P$ such that $PQ- \one =R_{-\infty}\in \cW^{-\infty}$.   $P^{-1}= Q- P^{-1}R_{-\infty}$ is invertible in $B(L^{2}(M; E))$, hence there exists $n\gg 1$ such that $Q'= Q- Q_{n}R_{-\infty}$ is invertible. We have $PQ'= \one+ R_{-\infty}- P Q_{n}R_{-\infty}= \one + R'_{-\infty}$. Since $P, Q_{n}\in \Psi^{0}$ and hence preserve Sobolev spaces, we obtain that $PQ_{n}R_{-\infty}\in \cW^{-\infty}$ hence $R'_{-\infty}\in \cW^{-\infty}$.  Since $\one + R'_{-\infty}= PQ'$, $\one + R'_{-\infty}$ is invertible in $B(L^{2}(M; E))$, hence by spectral invariance we have $(\one + R'_{-\infty})^{-1}\in \cW^{-\infty}$ and $P^{-1}= Q'(\one + R'_{-\infty})^{-1}\in \Psi^{0}$.  The statement about principal symbols is easy since $\sigma_{\rm pr}(P^{-1})= \sigma_{\rm pr}(Q')= (\sigma_{\rm pr}(P))^{-1}$. \qed
 
The same result is true in the time dependent case, with a similar proof.
  \begin{proposition}\label{prop4.1time}
Let $A(t)\in \cinfb(I;\Psi^{m}(M; L(E)))$, $m\geq 0$ elliptic. Then the following holds:
\ben
\item  $A(t)$ is closeable on $\coinf(M; E)$ with $\Dom A^{\rm cl}(t)= H^{m}(M; L(E))$.  
\item  if  there exists $\delta>0$ such that $[-\delta, \delta]\cap  \sigma(A^{\rm cl}(t))=\emptyset$ for $t\in I$, then $A^{-1}(t)\in \cinfb(I;\Psi^{-m}(M; L(E)))$ and
\[
\sigma_{\rm pr}(A^{-1})(t)= (\sigma_{\rm pr}(A))^{-1}(t).
\]
\een
\end{proposition}
\subsection{Functional calculus}
 \subsubsection{Elliptic selfadjoint operators}\label{sec4.2.5}
As in Subsect. \ref{sec4.2} we fix  a bounded Hilbertian structure $(\cdot| \cdot)_{E}$ on the fibers of $E$ and  define the scalar product
 \[
 (u|v)= \int_{M}(u(x)| v(x))_{E}dVol_{g}, \ u, v\in \coinf(M; E).
 \]
Let  $H(t)\in \cinfb(I;\Psi^{m}(M; L(E)))$ be  elliptic, symmetric on $\coinf(M; E)$. Using Prop. \ref{prop4.1} one easily shows that  its closure is selfadjoint with domain $H^{m}(M; E)$. Note also that its principal symbol $\sigma_{\rm pr}(H)(t, x, \xi)$ is selfadjoint for the Hilbertian scalar product on $E_{x}$. 

 We now discuss some results on the functional calculus for selfadjoint pseudodifferential operators.  As usual the functional calculus starts with the resolvent $(z-H)^{-1}$ and the properties of $z\mapsto (z- H)^{-1}$ can be obtained by considering pseudodifferential operators depending on some large parameter.
 
\subsubsection{Pseudodifferential operators with parameters}
 Let $(M, \hat{g})$ be of bounded geometry and $E\xrightarrow{\pi}M$ a vector bundle over $M$ of bounded geometry. We first discuss pseudodifferential operators on $M$ depending on some large parameter $\lambda\in \rr$.  
 
 We denote by $\widetilde{S}^{m}(T^{*}M; L(E))$ the space of symbols $b\in \cinf(\rr\times T^{*}M; L(E))$  such that if $b_{i}(\lambda)= T_{i}b(\lambda)$ are the  pushforwards of $b(\lambda)$ associated to a covering $\{U_{i}\}_{i\in \nn}$ (see \ref{sec4.1.1}), we have: 
  \[
\p^{\gamma}_{\lambda}\p^{\alpha}_{x}\p^{\beta}_{\xi}b_{i}(\lambda, x, \xi)\in O(\langle \xi\rangle + \langle \lambda\rangle)^{m -|\beta|- \gamma},\ (\lambda, x, \xi)\in \rr\times T^{*}B(0, 1)
  \]
  uniformly with respect to $i\in \nn$.   We denote by $\widetilde{S}^{m}_{\rm h}(T^{*}M; L(E))$ the subspace of such symbols which are homogeneous w.r.t. $(\lambda, \xi)$ and by $\widetilde{S}^{m}_{\rm ph}(T^{*}M; L(E))$ the subspace of polyhomogeneous symbols.

  We define  $\Op(b)u$ for $u\in \coinf(\rr_{\lambda}\times M; E)$ by:
  \[
  \Op(b)u(\lambda)\defeq  \Op(b(\lambda))u(\lambda).
  \]
  We define $\widetilde{\cW}^{-\infty}(M; L(E))$ as the set of smooth functions $b: \rr\in \lambda\mapsto b(\lambda)\in \cW^{-\infty}(M; L(E))$ such that
  \[
  \| \p^{\gamma}_{\lambda}b(\lambda)\|_{B(H^{-m}(\Sigma), H^{m}(\Sigma))}\in O(\langle \lambda\rangle^{-n}), \ \forall, m,n,\gamma\in \nn,
  \]
  and set
  \[
  \widetilde{\Psi}^{m}(M; L(E))\defeq  \Op (\widetilde{S}_{\rm ph}^{m}(T^{*}M; L(E)))+ \widetilde{\cW}^{-\infty}(M; L(E)).
  \]
  
  We also define the time dependent versions of the above spaces:
  \[
  \cinfb(I; \widetilde{S}^{m}_{\rm ph}(T^{*}M; L(E))), \ \cinfb(I; \widetilde{\cW}^{-\infty}(M; L(E))), \cinfb(I; \widetilde{\Psi}^{m}(M; L(E))).
  \]
  We define the {\em principal symbol} of  $A(t)\in \cinfb(I; \widetilde{\Psi}^{m}(M; L(E)))$ as in \ref{sec4.2.3},  using the polyhomogeneity.
  
  An operator $A(t)\in \cinfb(I; \widetilde{\Psi}^{m}(M; L(E)))$ is {\em elliptic} if  $\sigma_{\rm pr}(A)(t)$ is invertible for $t\in I$ and $\sigma_{\rm pr}(A)^{-1}(t)\in \cinfb(I; \widetilde{S}^{-m}_{\rm ph}(M; L(E)))$.

\begin{proposition}\label{propopi}
 Let $H(t)\in \cinfb(I; \Psi^{1}(M; L(E)))$ elliptic and formally selfadjoint.  Let us still denote by $H(t)$ its closure, which is selfadjoint on $H^{1}(M; E)$ by  Prop. \ref{prop4.1time}. Assume that there exists $\delta>0$ such that $[-\delta, \delta]\cap \sigma(H(t))= \emptyset$ for $t\in I$.
 
  Then $\rr\ni \lambda\mapsto (H(t)+ \i \lambda)^{-1}$ belongs to $\cinfb(I; \widetilde{\Psi}^{-1}(M; L(E)))$ with principal symbol $(\sigma_{\rm pr}(H(t))+ \i \lambda)^{-1}$.
\end{proposition}
  \proof We will use   arguments similar to those used  in \cite[Appendix A.1]{GW4}, reducing ourselves to 
 the situation without parameters.  
 
  {\em Step 1.} we   introduce a variable $l$ dual to $\lambda$ and 
 equip $\widetilde{M}=\rr_{l}\times M$ with  the metric $\widetilde{g}= dl^{2}+ \hat{g}$.  $(\widetilde{M}, \widetilde{g})$ is of bounded geometry, so we can define the spaces $S^{m}_{\rm ph}(T^{*}\widetilde{M}; L(E))$. 
   
   As  bounded atlas of $\widetilde{M}$ we can take $\widetilde{U}_{i}= \rr\times U_{i}$, $\widetilde{\psi}_{i}(l, x)= (l, \psi_{i}(x))$ for $(U_{i}, \psi_{i})_{i\in\nn}$ a bounded atlas of $M$, and as bounded partition of unity on $\widetilde{M}$ a bounded partition of unity on $M$.

   This implies that the  quantization map $\widetilde{\Op}$  obtained from such choices is the usual quantization on $\rr$ in the variables $(l, s)$, tensorized with the quantization map $\Op$ on $M$.
   
  The appropriate ideal of smoothing operators, dictated by the definition of $\widetilde{\cW}^{-\infty}(M; L(E))$ is now
 \[
\widetilde{\cW}^{-\infty}(\widetilde{M}; L(E))= \{A: {\rm ad}^{n}_{l}A\in B(H^{-m}(\widetilde{M}), H^{m}(\widetilde{M})), \ \forall m,n\in \nn\}.
 \]
 Note that this ideal is smaller than the one defined in \ref{sec4.2.1b} (with  $M$ replaced by $\widetilde{M}$) because of the additional control on the multicommutators ${\rm ad}^{n}_{l}A$.

 We define then
 \[
 \widetilde{\Psi}^{m}(\widetilde{M}, L(E))= \widetilde{\Op}(S^{m}_{\rm ph}(T^{*}\widetilde{M}; L(E)))+ \widetilde{\cW}^{-\infty}(\widetilde{M}; L(E)).
 \]
 As is now usual we introduce also the time dependent versions $ \cinfb(I;\widetilde{\Psi}^{m}(\widetilde{M}, L(E)))$ etc.
 
{\em Step 2.} we study the link between $ \cinfb(I;\widetilde{\Psi}^{m}(\widetilde{M}; L(E)))$ and $ \cinfb(I;\widetilde{\Psi}^{m}(M; L(E)))$. We note first that
 \[
  \cinfb(I;\widetilde{S}^{m}_{\rm ph}(T^{*}M; L(E)))= \{b(t)\in  \cinfb(I;S^{m}_{\rm ph}(T^{*}\widetilde{M}; L(E))): \p_{l}b(t)=0\},
 \]
 and denoting by $T_{l}$ the group of translations  in $l$, we have:
 \[
 [T_{l}, \widetilde{\Op}(b)(t)]=0\ \forall l\in \rr\Leftrightarrow b(t)\in  \cinfb(I;\widetilde{S}^{m}_{\rm ph}(T^{*}M; L(E))).
 \]
 Equivalently, if $\cF: L^{2}(\rr, dl)\tosim L^{2}(\rr, d\lambda)$ is the unitary Fourier transform in $l$, we have
\begin{equation}
\label{a2}
\begin{array}{rl}
&b(t)\in  \cinfb(I;S^{m}_{\rm ph}(T^{*}\widetilde{M}; L(E))), \ [T_{l}, \widetilde{\Op}(b)(t)]=0\ \forall l\in \rr\\[2mm]
\Leftrightarrow& \cF \widetilde{\Op}(b)(t)\cF^{-1}= \int_{\rr}^{\oplus}\Op (b(t, \lambda))dl, \hbox{ for }b(t, \cdot)\in  \cinfb(I;\widetilde{S}^{m}(T^{*}M; L(E))).
\end{array}
\end{equation}
As in \cite[Appendix A.1]{GW4} we prove similarly that 
\begin{equation}
\label{a4}
\begin{array}{rl}
&w(t)\in  \cinfb(I;\widetilde{\cW}^{-\infty}(\widetilde{M}; L(E))), \ [w(t), T_{\lambda}]=0, \ \forall \lambda\in \rr\\[2mm]
\Leftrightarrow& \cF w(t)\cF^{-1}= \int_{\rr}^{\oplus}w(t,\lambda)d\lambda\hbox{ for }w(t, \cdot)\in  \cinfb(I;\widetilde{\cW}(M; E)).
\end{array}
\end{equation}
{\em Step 3.} conjuguating by $\cF$, we consider the operator $A(t)= H(t)+ \i D_{l}$, which is elliptic in $ \cinfb(I;\Psi^{1}(\widetilde{M}; L(E)))$. Its closure is a normal operator with domain $H^{1}(\widetilde{M}; E)$ and by assumption there exists $\delta>0$ such that $[-\delta, \delta]\cap \sigma(A(t))= \emptyset$.

 We check the spectral invariance of  the ideal $ \cinfb(I;\widetilde{\cW}^{-\infty}(\widetilde{M}; L(E)))$ as in \cite[Appendix A.1]{GW4}. We can hence  apply the abstract results in \cite{alnv}.  We obtain that $A(t)^{-1}\in \cinfb(I; \Psi^{-1}(\widetilde{M}; L(E)))$ with principal symbol $(\sigma_{\rm pr}(H(t))+ \i \lambda)^{-1}$. Moreover since $[A(t), T_{l}]=0$ we also obtain that $[A(t)^{-1}, T_{l}]=0$.  We have
 \[
 \cF A(t)\cF^{-1}= \int_{\rr}^{\oplus}(H(t)+ \i \lambda)d\lambda,
 \]
 and hence
 \[
 \cF A(t)^{-1}\cF^{-1}= \int_{\rr}^{\oplus} (H(t)+ \i \lambda)^{-1}d\lambda.
 \]
 By Step 2 we obtain that  $\lambda\mapsto (H(t)+ \i \lambda)^{-1}$ belongs to $\cinfb(I; \widetilde{\Psi}^{-1}(M; L(E)))$, with principal symbol 
  $(\sigma_{\rm pr}(H(t))+ \i \lambda)^{-1}$. \qed

Let us now prove a useful consequence of Prop. \ref{propopi}.
 \begin{proposition}\label{prop4.3}
Let  $H(t)\in \cinfb(I; \Psi^{1}(M; L(E)))$ be elliptic, symmetric on $\coinf(M; E)$, and let us denote still by $H(t)$ its closure.  Assume that there exists $\delta>0$ such that $[-\delta, \delta]\cap \sigma(H(t))= \emptyset$ for $t\in I$. Then the spectral projections $\one_{\rr^{\pm}}(H(t))$ belong to $\cinfb(I;\Psi^{0}(M; L(E)))$ and
\[
\sigma_{\rm pr}(\one_{\rr^{\pm}}(H(t)))= \one_{\rr^{\pm}}(\sigma_{\rm pr}(H(t))).
\]
\end{proposition}
\proof 
Since $(H^{2}+ \lambda^{2})^{-1}= (H+\i \lambda)^{-1}(H -\i \lambda)^{-1}$, we obtain from Prop. \ref{propopi} that $\lambda\mapsto (H^{2}(t)+ \lambda^{2})^{-1}$ belongs to $\cinfb(I; \tilde{\Psi}^{-2}(M; L(E)))$, with principal symbol $(\sigma_{\rm pr}H(t)^{2}+ \lambda^{2})^{-1}$. Using
\[
a^{-\12}= \frac{2}{\pi}\int_{0}^{+\infty}(a+s^{2})^{-1}ds, \ a>0
\]
we obtain that $|H(t)|^{-1}\in\cinfb(I; \Psi^{-1}(M; L(E)))$ with principal symbol $|\sigma_{\rm pr}(H)(t)|^{-1}$. We write then 
$\one_{\rr^{\pm}}(H)= \12(\one \pm \frac{H}{|H|})$ to obtain the proposition. \qed

%
%

  \section{Pure Hadamard states for Dirac fields}\label{sec5}\init
  In this section we prove the main result of this paper, namely the construction of pure Hadamard states for Dirac fields on globally hyperbolic spacetimes of bounded geometry.  We will work in the framework described in Subsect. \ref{sec5.0}. 
\subsection{Reduction of Dirac operators on cartesian products}\label{sec5.1}
Let $D$ be a Dirac operator on the spacetime $(M, g)$ considered in \ref{carproduct}, i.e. $M= I\times \Sigma$ and $g= - dt^{2}+ h_{t}(\rx)d\rx^{2}$.

We will denote by $(\rx, k)$ local coordinates on $T^{*}\Sigma$, and by $(x, \xi)= ((t, \rx), (\tau, k))$ local coordinates on $T^{*}M$.

We recall that   the restriction  $S_{t}(\Sigma)$  of the spinor bundle $S(M)$ to $\Sigma_{t}$ is independent on $t$ and denoted by $S(\Sigma)$.  

We recall also that $S(M)$ is equipped by a Hermitian form $\beta$,  see 
\ref{sec1.3.3} and we denote by $\beta_{t}$ its  restriction  to $S(\Sigma)$. Also  we denote by $\gamma_{t}: T_{\Sigma_{t}}M\to L(S(\Sigma))$ the restrictions of $\gamma$ to $S(\Sigma)$.

 We first simplify the Dirac operator $D$ using  parallel transport by  $\p_{t}$ for the spin connection.
 
  For $f\in \cinf(\Sigma_{s}; S(\Sigma))$ we denote by $\mathcal{T}(s)f= \psi$ the solution of
\beq\label{e3.8}\left\{
\begin{array}{l}
\nabla^{S}_{\p_{t}}\psi=0\hbox{ in }I\times \Sigma\\
\psi_{| \Sigma_{s}}=f,
\end{array}\right.
\eeq
and set
\beq\label{e3.8b}
\begin{array}{l}
\mathcal{T}(t, s)f= \mathcal{T}(s)f_{| \Sigma_{t}},\\[2mm]
\mathcal{T}: \cinf(I; \cinf(\Sigma, S(\Sigma)))\to \cinf(M; S(M))\\
\psi(t)\mapsto (\mathcal{T}\psi)(t)= |h_{t}|^{-\frac{1}{4}}|h_{0}|^{\frac{1}{4}}\mathcal{T}(t, 0)\psi(t),
\end{array}
\eeq
\begin{lemma}\label{lemma3.1}
 One has
 \ben
 \item $\cT(s,t)\gamma_{t}(e_{0})\cT(t, s)= \gamma_{s}(e_{0})$, $t,s\in I$.
\item $\cT(s,t)\gamma_{t}(e_{a}(t))\cT(t, s)= \gamma_{s}(e_{a}(s))$, $t,s\in I$.
\item $\cT(t,s)^{*}\beta_{t}\cT(t,s)= \beta_{s}$, $t,s\in I$.
 \een
\end{lemma}
\proof Since $\nabla_{e_{0}}e_{a}=0$, we obtain using \eqref{e10.4} {\it i)} that
\[
\nabla_{e_{0}}^{S} \gamma(e_{a})\cT(s)f= \gamma(e_{a})\nabla_{e_{0}}^{S}\cT(s)f=0,
\]
hence 
\[\cT(s,t)\gamma_{t}(e_{0})\cT(t, s)= \gamma_{s}(e_{0}), \ \cT(s,t)\gamma_{t}(e_{a}(t))\cT(t, s)= \gamma_{s}(e_{a}(s))
\]  as claimed. Similarly using  \eqref{e10.4} {\it ii)} we have
\[
\p_{t}(\bar{\cT(t,s )f})\dual \beta_{t}\cT(t,s)f= \bar{\nabla^{S}_{e_{0}}\cT(s)f}\dual \cT(s)f+ \bar{\cT(t,s )f})\dual \beta_{t}\nabla^{S}_{e_{0}}\cT(t,s)f=0,
\]
which proves (3). \qed
\begin{proposition}\label{prop3.1}
Let
\[
\cD\defeq \cT^{-1}(\gamma(e_{0})D)\cT,
\]
and
\[
\bar{\tilde{\psi}}\dual \tilde{\nu}\tilde{\psi}\defeq \i \int_{\Sigma}\bar{\tilde{\psi}}\dual \beta_{0}\gamma_{0}(e_{0})\tilde{\psi}|h_{0}|^{\12}d\rx
\]
Then
\ben
\item the map
 \[
 \cT: (\Sol (\cD), \tilde{\nu})\tosim (\Sol(D), \nu)
 \]
 is unitary.
 \item  We have 
\[
\cD= \p_{t}- \i H(t),
\]
where $H(t)\in\cinf_{\rm b}(I; \Psi^{1}(\Sigma, S(\Sigma)))$,
\beq\label{e3.1}
\sigma_{\rm pr}(H(t))(\rx, k)= - \gamma_{0}(e_{0})\gamma_{0}(h_{t}(\rx)^{-1}k)
\eeq
and $H(t)$ is formally selfadjoint for the scalar product
\[
\bar{f}\dual \tilde{\nu} f\defeq \i\int_{\Sigma}\bar{f}\dual \beta_{0}\gamma_{0}(e_{0})f |h_{0}|^{\12}d\rx.
\]
\een
\end{proposition}
\proof
(1) is obvious since $\cT(0, 0)=\one$.  Let us now prove (2).
From \eqref{e3.8} we obtain 
\begin{equation}
\label{e3.9}
\mathcal{T}\p_{t}\mathcal{T}^{-1}=\nabla_{e_{0}}^{S}- \frac{1}{4}\p_{t}|h_{t}||h_{t}|^{-1}.
\end{equation}
If we fix a local oriented and time oriented orthonormal frame $(e_{a})_{0\leq a \leq d}$ over some open set $U= I\times V$, we have
\[
D= - \gamma(e_{0})\nabla^{S}_{e_{0}}+ \gamma(e_{a})\nabla^{S}_{e_{a}}+m,
\]
(where we sum only for $1\leq a\leq d$),
and
\beq\label{e1.1bd}
\begin{array}{l}
\mathcal{T}^{-1}\gamma(e_{0})D \mathcal{T}= \p_{t}- \i H(t), \\[2mm]
 H(t)= \i \cT^{-1}(\gamma_{t}(e_{0})\gamma_{t}(e_{a}(t))\nabla_{e_{a}(t)}^{S}+ \gamma_{t}(e_{0})m)\mathcal{T}+ \frac{1}{4}\p_{t}|h_{t}||h_{t}|^{-1}.
\end{array}
\eeq
Let us now prove the properties of $H(t)$ stated in (2).
Let us first consider the  operator 
\[
\tilde{H}(t)=  \i \gamma_{t}(e_{0})\gamma_{t}(e_{a}(t))\nabla_{e_{a}(t)}^{S}+ \i  \gamma_{t}(e_{0})m
\]
appearing in \eqref{e1.1bd}. 

We fix a bounded atlas $(V_{i}, \psi_{i})_{i\in \nn}$ of $(\Sigma, h_{0})$ and set $U_{i}= I\times V_{i}$.   We use the local frames $\cE_{i}= (e_{i, a})_{0\leq a\leq d}$ as in Subsect. \ref{spinprod},
and  the spin frames $\cB_{i}= (E_{i, A})_{1\leq A\leq N}$ associated to the frames $\cE_{i}$ over $U_{i}$. The family  $(\cE_{i})_{i\in \nn}$ resp.  $(\cB_{i})_{i\in \nn}$ is a bounded family of frames of $TM$, resp. $S(M)$ over $U_{i}$. 

An easy computation shows that $\Gamma^{c}_{ab}= \nabla_{e_{a}}e_{b}\dual e^{c}$ belongs to $\cinfb(I\times V)$,
and  if we reintroduce the index $i$ and set $V= V_{i}$, then the seminorms of $\Gamma^{c}_{ab}$  in $\cinfb(I\times V)$ are uniform with respect to $i\in \nn$.  Using the expression \eqref{e2.7} of $\nabla^{S}$ in the spin frames $\cB_{i}$ we obtain that $\tilde{H}(t)\in \cinfb(I; {\rm Diff}^{1}(\Sigma_{t}; L(S(\Sigma))))$.

Let us now consider the maps $\cT(t,s)$.
Let us  again forget the index $i$ and  denote by $\pmb{\psi}\in \cc^{N}$ the components of $\psi$ in the frame $\cB$. Then from \eqref{e2.7}, we have:
\[
\nabla^{S}_{e_{0}}\pmb{\psi}= \p_{t}\pmb{\psi}+ \frac{1}{4} \Gamma^{a}_{0b} \pmb{\gamma}_{a}\pmb{\gamma}^{b}\pmb{\psi}.
\]
Using elementary bounds on solutions of linear differential equations, we conclude that $I^{2}\ni (t,s)\mapsto \cT(t,s)\in L(S(\Sigma_{s}), S(\Sigma))$ belongs to $\cinfb(I^{2}; L(S(\Sigma), S(\Sigma)))$.

This implies  that $H(t)\in \cinf_{\rm b}(I; \Psi^{1}(\Sigma, S(\Sigma)))$. Its principal symbol  is clearly given by  \eqref{e3.1}. 

We know that if $D\psi= 0$ then
\[
\int_{\Sigma}\bar{\psi}(t)\dual \beta_{0} \gamma_{0}(e_{0}) \psi(t)|h_{t}|^{\12}d\rx
\]
is independent on $t$, hence if $\tilde{\psi}= \mathcal{T}^{-1}\psi$ we obtain using Lemma \ref{lemma3.1} that
\[
\int_{\Sigma}\bar{\tilde{\psi}}(t)\dual \beta_{0} \gamma_{0}(e_{0}) \tilde{\psi}(t)|h_{0}|^{\12}d\rx
\]
is independent on $t$. This implies that $H(t)= H^{*}(t)$ on $\coinf(\Sigma; S(\Sigma))$  for the scalar product $\tilde{\nu}$. \qed

\subsection{Model case}\label{sec5.2}
 In the previous section we have reduced the Dirac equation to the equation:
 \[
 \p_{t}\psi- \i H(t)\psi=0,
 \]
for   $H(t)$  as in Prop. \ref{prop3.1}.  The space $\coinf(I\times \Sigma; S(\Sigma))$ is equipped with the Hilbertian scalar product
\[
\bar{\psi}\dual \nu \psi= \int_{I\times \Sigma}\bar{\psi}\dual \beta_{0}\psi dt|h_{0}|^{\12}dx,
\]
while $\coinf(\Sigma; S(\Sigma))$ is equipped with
\[
\bar{f}\dual \nu f= \int_{\Sigma}\bar{f}\dual \beta_{0}f |h_{0}|^{\12}dx.
\]
Adjoints of operators will always be computed with respect to these scalar products. Our reference Hilbert space is
\[
\cH= L^{2}(\Sigma; S(\Sigma)),
\]
equal to the completion of $\coinf(\Sigma; S(\Sigma))$  for $\nu$.  As in \ref{sec4.2.5} $H(t)$ with domain $H^{1}(\Sigma, S(\Sigma))$ is selfadjoint on $\cH$.

We will need later to consider the spectral projections $\one_{\rr^{\pm}}(H(t))$, which can be singular  if $0\in \sigma(H(t))$. In the next lemma we show how to modify $H(t)$ by a time-dependent smoothing operator to avoid this problem.
\begin{lemma}\label{lemma5.1}
 There exists $R_{-\infty}\in \cinfb(I; \cW^{-\infty}(\Sigma; S(\Sigma)))$ with $R_{-\infty}(t)= R_{-\infty}(t)^{*}$ such that
 \[
 \sigma(H(t)+ R_{-\infty}(t))\cap [-1, 1]= \emptyset, \ t\in I.
 \]
 \end{lemma}
\proof 
Let us forget the parameter $t$ to ease notation. Using the principal symbol of $H(t)$ and the Clifford relations we obtain that
$ H^{*}H= h_{2}+ R_{1}$, 
where $h_{2}\in \Psi^{2}(\Sigma, S(\Sigma))$ is a scalar selfadjoint  operator with principal symbol $k\dual h_{t}(x)k$ and $R_{1}\in \Psi^{1}(\Sigma; S(\Sigma))$. 

 Let $\chi\in \coinf(\rr)$ a cutoff function equal to $1$ near $0$ and
$\chi_{\lambda}= \chi(\lambda^{-2}h_{2})$.
If  $H_{\lambda}= H + \i \lambda \chi_{\lambda} \gamma_{0}$, then $H_{\lambda}^{*}= H_{\lambda} $ and
\[
\begin{array}{rl}
H_{\lambda}^{*}H_{\lambda}=& H^{*}H+ \lambda^{2}\chi_{\lambda}^{2}+ \i \lambda (\chi_{\lambda} \gamma_{0}H + H \chi_{\lambda} \gamma_{0})\\[2mm]
=& H^{*}H+ \lambda^{2}\chi_{\lambda}^{2}+ \lambda[H, \chi_{\lambda}]\gamma_{0}+ \chi_{\lambda} (H\gamma_{0}+ \gamma_{0}H).
\end{array} 
\]
Since $\chi_{\lambda}$ is scalar, $[H, \chi_{\lambda}]\in \Psi^{0}$ with  $\| [H, \chi_{\lambda}]\| \in O(1)$ uniformly for $\lambda \gg 1$ and using the Clifford relations $H\gamma_{0}+ \gamma_{0}H\in \Psi^{0}$.  Therefore we have
\[
H_{\lambda}^{*}H_{\lambda}= h_{2}+ \lambda ^{2}\chi(\lambda^{-2}h_{2})+ R_{1}+ \lambda R_{0}(\lambda),
\]
where $R_{1}\in \Psi^{1}$ and $\| R_{0}(\lambda)\| \in O(1)$. For $\lambda$ large enough we have $h_{2}+  \lambda^{2}\chi^{2}(\lambda^{-2}h_{2})\geq \12 (h_{2}+ \lambda^{2})$. Since $R_{1}\geq - \epsilon h_{2}- C\epsilon^{-1}$ for all $\epsilon>0$ can pick $\lambda$ large enough so that $H_{\lambda}^{*}H_{\lambda}\geq 2$ and hence $[-1, 1]\cap \sigma(H_{\lambda})= \emptyset$. To finish the proof we note that $\chi_{\lambda}\gamma_{0}\in \cW^{-\infty}$. \qed

\subsubsection{Unitary group}
Let us denote by $\cU(t,s)$, $s,t\in I$ the unitary evolution generated by $H(t)$, i.e. the solution of
\[
\left\{
\begin{array}{l}
\p_{t}\cU(t, s)= \i H(t)\cU(t, s), \\
 \p_{s}\cU(t,s)= - \i \cU(t, s)H(s),\\
\cU(s, s)= \one.
\end{array}\right.
\]
The properties of $H(t)$ imply that $\cU(t,s)$ is well-defined by a classical result of Kato, see for example \cite{SG}. For later use, let us first prove an easy result.
\begin{lemma}\label{lemma5.2}
 $\cU(t,s)$ are uniformly bounded in $B(H^{m}(\Sigma; S(\Sigma)))$ for $t,s\in I$, $m\in \rr$.
\end{lemma}
\proof
Let  $\epsilon= (-\Delta_{h_{0}}+ 1)^{\12}$. For $u\in \epsilon^{m}\coinf(\Sigma; S(\Sigma))$ we set
$f(t)= \|\cU(s,t)\epsilon^{m}\cU(t, s)\epsilon^{-m}u\|$, which is finite since $\cU(t,s)$ preserves $\coinf(\Sigma; S(\Sigma))$.  We have
\[
\begin{array}{rl}
|f'(t)|\leq& \|\cU(s,t)[H(t), \i \epsilon^{m}]\cU(t,s)\epsilon^{-m}u\| \\[2mm]
\leq& \| \cU(s,t)[H(t), \i \epsilon^{m}]\epsilon^{-m}\cU(t, s)\|f(t)\leq C f(t), \ t\in I.
\end{array}
\]
By Gronwall's inequality we obtain that $f(t)\leq C f(s)$ for $t,s\in I$ hence 
\[
\|\cU(s,t)\epsilon^{m}\cU(t, s)\epsilon^{-m}u\| \leq C \| u\|,\ u\in \epsilon^{m}\coinf(\Sigma; S(\Sigma)),
\]
 which proves the lemma since $\epsilon^{m}\coinf(\Sigma; S(\Sigma))$ is dense in $L^{2}(\Sigma; S(\Sigma))$. \qed

\subsection{Some preparations}\label{sec5.4}
 We start by describing a concise framework to solve  recursive equations that are often encountered in symbolic calculus. A similar method was used in \cite[Lemma A.1]{GW1}.
\begin{definition}\label{def5.1}
 For $p\in \nn$ we denote by $\cF_{-p}$ the set of maps
 \[
 F: \cinfb(I; \Psi^{-1}(\Sigma; S(\Sigma)))\to \cinfb(I; \Psi^{-p}(\Sigma; S(\Sigma)))\hbox{ such that}
 \]
 $F(R_{1})- F(R_{2})\in \cinfb(I; \Psi^{-p-j}(\Sigma; S(\Sigma)))$ if $R_{1}- R_{2}\in \cinfb(I; \Psi^{-1-j}(\Sigma; S(\Sigma)))$.
 \end{definition}
One can call an element of $\cF_{-p}$ a {\em symbolic contraction} of order $p$.  The following proposition is proved exactly as \cite[Lemma A.1]{GW1}.
\begin{proposition}\label{prop5.1}
 Let $A\in \cinfb(I; \Psi^{-1}(\Sigma; S(\Sigma)))$ and $F_{-2}\in \cF_{-2}$. Then there exists a solution $R\in \cinfb(I; \Psi^{-1}(\Sigma; S(\Sigma)))$, unique modulo $\cinfb(I; \cW^{-\infty}(\Sigma; S(\Sigma)))$ of the equation:
 \[
 R= A+ F_{-2}(R)\hbox{ mod }\cinfb(I; \cW^{-\infty}(\Sigma; S(\Sigma))).
 \]
 \end{proposition}
\proof Let us denote $F_{-2}$ simply by $F$ and $\cinfb(I; \Psi^{p}(\Sigma; S(\Sigma)))$ simply by $\Psi^{p}$. We set $S_{0}= A$, $S_{n}= A+ F(S_{n-1})$ for $n\geq 1$. We have $S_{1}- S_{0}= F(A)$ and $S_{n}- S_{n-1}= F(S_{n-1})- F(S_{n-2})$. Since $F\in \cF_{ -2}$ we obtain by induction that $S_{n}- S_{n-1}\in \Psi^{ -(n+1)}$. We take $R\in \Psi^{-1}$ such that $R\sim S_{0}+ \sum_{0}^{\infty}S_{n}- S_{n-1}$ which solves the equation modulo $\Psi^{ -\infty}$. If $R_{1}, R_{2}$ are two solutions then
 $R_{1}- F_{2}= F(F_{1})- F(R_{2})$ modulo $\Psi^{ -\infty}$ hence using that $F\in \cF_{ -2}$ we obtain by induction on $n$ that $R_{1}- R_{2}\in \Psi^{ -n}$ for all $n\in \nn$ which proves uniqueness modulo $\Psi^{ -\infty}$. \qed

 An element of $\cF_{-p}$ will be denoted by $F_{-p}$. Let us collect  some easy properties.
 
\begin{lemma}\label{lemma5.2b}
 \ben
 \item If $A\in \cinfb(I; \Psi^{k}(\Sigma; S(\Sigma)))$  and $F_{-p}\in \cF_{-p}$ then  the maps \[
 \begin{array}{l}
 A F_{-p}: R\mapsto AF_{-p}(R)\\[2mm]
 F_{-p}A: R\mapsto F_{-p}(R)A
 \end{array}
 \]
  belong to $\cF_{-p+k}$ for $k\leq p$.
 \item If $F_{-p_{i}}\in \cF_{-p_{i}}$ then  the map \[
 F_{-p_{1}}F_{-p_{2}}: R\mapsto F_{-p_{1}}(R)F_{-p_{2}} (R)\]
  belongs to $\cF_{-p_{1}- p_{2}}$.
  \item the map $R\mapsto R^{p}$ belongs to $\cF_{-p}$ for $p\in \nn$.
 \item the map $R\mapsto \e^{R}$ belongs to $\cF_{0}$.
\item one has $\e^{R}= 1 + R + F_{-2}(R)$, where $F_{-2}\in \cF_{-2}$.
 \een
\end{lemma}
\proof (1) and (2) are easy. Since the identity map belongs to $\cF_{-1}$, (3) follows from (2). If $R\in \cinfb(I; \Psi^{-1})$ then $\e^{R}= \sum_{n\in \nn}\frac{1}{n!}R^{n}$ as a norm convergent series. Since $R^{n}\in \cinfb(I; \Psi^{-n})$,  we obtain that $\e^{R}\in \cinfb(I; \Psi^{-0})$.
Next we write
\[
\e^{ R_{1}}- \e^{R_{2}}=  \int_{0}^{1}\e^{ \theta R_{1}}(R_{1}- R_{2})\e^{(1- \theta)R_{2}}d\theta.
\]
This implies that $\e^{ R_{1}}- \e^{R_{2}}\in \cinfb(I; \Psi^{-j})$ if $R_{1}- R_{2}\in \cinfb(I; \Psi^{-j})$ and proves (4). We have also
\[
\e^{R}= \one+ R+ \int_{0}^{1}(1- \theta)R^{2}\e^{\theta R}d\theta.
\]
By (3) the map $R\mapsto R^{2}$ belongs to $\cF_{-2}$ and by (4) the map $R\mapsto \e^{ \theta R}$ belongs to $\cF_{-0}$. This implies (5), using (2). \qed

\subsection{Construction of some projections}\label{sec5.3}

\begin{proposition}\label{prop5.2}
 There exists time-dependent projections $\tilde{P}^{\pm}\in \cinfb(I; \Psi^{0}(\Sigma, L(S(\Sigma))))$ such that:
\ben
\item $\tilde{P}^{\pm}(t)= \tilde{P}^{\pm}(t)^{*}$, $\tilde{P}^{+}(t)+ \tilde{P}^{-}(t)=\one$;
\item $\cU(t, s)\tilde{P}^{\pm}(s)= \tilde{P}^{\pm}(t)\cU(t,s)+ R_{-\infty}(t,s)$, 

where $R_{-\infty}(t,s)\in \cinfb(I^{2}; \cW^{-\infty}(\Sigma; L(S(\Sigma))))$;
\item $\WF(\cU(\cdot, s)\tilde{P}^{\pm}(s))'\subset (\cN^{\pm}\cup \cF)\times T^{*}\Sigma$, for $\cF= \{k=0\}\subset T^{*}M$.
\een 
\end{proposition}
\proof
The construction of $\tilde{P}^{\pm}$ is divided in several steps.

{\it Step 1}. In the first step we  replace $H(t)$ by $\tilde{H}(t)=H(t)+ R_{-\infty}(t)$ so that $[-1, 1]\cap \sigma(\tilde{H}(t))= \emptyset$, using Lemma \ref{lemma5.1}. Let $ \tilde{\cU}(t,s)$ the unitary group with generator $\tilde{H}(t)$. From Lemma \ref{lemma5.2} and Duhamel's formula we obtain easily that $ \tilde{\cU}(t,s)- \cU(t,s)\in  \cinfb(I^{2}; \cW^{-\infty}(\Sigma; L(S(\Sigma))))$, so denoting $\tilde{H}(t)$ again by $H(t)$ we can assume without loss of generality that $[-\delta, \delta]\cap \sigma(H(t))= \emptyset$ for $t\in \rr$.

By Prop. \ref{prop4.3} the projections
\[
P^{\pm}(t)= \one_{\rr^{\pm}}(H(t))
\]
are well defined, selfadjoint with $P^{\pm}(t)\in \cinfb(I; \Psi^{0}(\Sigma; L(S(\Sigma))))$ and
\begin{equation}
\label{e5.0}
\sigma_{\rm pr}(P^{\pm})(t, \rx, k)= \one_{\rr^{\pm}}(\sigma_{\rm pr}(H)(t, \rx, k)).
\end{equation}
Since $\sigma_{\rm pr}(H(t, \rx, k))= -\gamma_{0}\gamma(h_{t}^{-1}(\rx)k)$, we obtain using the Clifford relations that:
\[
\sigma_{\rm pr}(P^{\pm})(t, \rx, k)\sigma_{\rm pr}(H(t, \rx, k))= \pm\epsilon(t, \rx, k)\sigma_{\rm pr}(P^{\pm})(t, \rx, k),
\]
for $\epsilon(t, x, k)= (k\dual h_{t}^{-1}(\rx)k)^{\12}$. By symbolic calculus this implies that
\begin{equation}
\label{e5.-0}
P^{\pm}(t)H(t)= \pm \epsilon(t, \rx, D_{x}) P^{\pm}(t)+ R_{0}^{\pm}(t), 
\end{equation}
where $R_{0}^{\pm}(t)\in \cinfb(I; \Psi^{0}(\Sigma; L(S(\Sigma))))$.

 Of course $\cU(t, s)P^{\pm}(s)\neq P^{\pm}(t)\cU(t,s)$ but we can try to modify $P^{\pm}(t)$ so that the equality holds, at least up to a smooth error term.

 {\it Step 2}. In the second step, we modify  $P^{\pm}(t)$.
 To this end let  $ R(t)\in \cinfb(I; \Psi^{-1}(\Sigma; L(S(\Sigma))))$ with $R(t)= R^{*}(t)$ and let us set
 \[
 \tilde{\cU}(t, s)\defeq \e^{\i R(t)}\cU(t,s)\e^{- \i R(s)}.
 \]
This  is a strongly continuous unitary group with generator:
 \[
 \tilde{H}(t)= \e^{\i R(t)}H(t)\e^{- \i R(t)}+ \i^{-1}\p_{t}\e^{\i R(t)}\e^{- \i R(t)}.
 \]
 The equation
 \[
 \tilde{\cU}(t,s)P^{\pm}(s)= P^{\pm}(t)\tilde{\cU}(t,s)
 \]
 is equivalent to
 \beq\label{e5.1}
 \p_{t}P^{\pm}(t)+ [P^{\pm}(t), \i \tilde{H}(t)]=0,
 \eeq
 and implies that if
 \beq\label{e5.2}
 \tilde{P}^{\pm}(t)= \e^{- \i R(t)}P^{\pm}(t)\e^{\i R(t)}
 \eeq
 one has
 \beq\label{e5.3}
  \tilde{P}^{\pm}(t)= \tilde{P}^{\pm}(t)^{*}, \ \cU(t, s)\tilde{P}^{\pm}(s)=\tilde{P}^{\pm}(t)\cU(t,s), \ t, s\in I.
 \eeq
 
{\it Step 3}. In the third step we solve equation \eqref{e5.1} for $R$, modulo a smoothing error.
For ease of notation we denote simply by $A$ a time-dependent pseudodifferential operator $A(t)$. 

Let $R\in \cinfb(I; \Psi^{-1})$.
By Lemma \ref{lemma5.2b} (5) we have $\p_{t}\e^{\i R}= \i \p_{t}R+ F_{-2}(R)$, $\e^{- \i R}= \one + F_{-1}(R)$ hence
\[
\i^{-1}\p_{t}(\e^{\i R})\e^{- \i R}= \p_{t}R+ F_{-2}(R)= F_{-1}(R).
\]
Using again Lemma \ref{lemma5.2b} (5) we also obtain
\[
\e^{\i R}H\e^{- \i R}= H + [R, \i H]+ F_{-1}(R)
\]
and hence
\beq\label{e5.3b}
\tilde{H}= H + [R, \i H]+ F_{-1}(R).
\eeq
We will look for $R$ of the form
\beq
\label{e5.4}
R=T(S)=  P^{+}SP^{-}+ P^{-}S^{*}P^{+}, \ S\in \cinfb(I; \Psi^{-1}).
\end{equation}
Note that if $F_{-p}\in \cF_{-p}$ then the map $S\mapsto F_{-p}(T(S))$ belongs also to $\cF_{-p}$ (note that the map $S\mapsto S^{*}$ belongs to  $\cF_{-1}$).

 We can now solve the equation \eqref{e5.1} modulo a smoothing error. Since $P^{\pm}$ are projections we have
 \[
 \begin{array}{rl}
& \p_{t}P^{\pm}+ [P^{\pm}, \i \tilde{H}]\\[2mm]
=& P^{+}(\p_{t}P^{\pm}+ [P^{\pm},  \i\tilde{H}])P^{-}+ P^{-}(\p_{t}P^{\pm}+ [P^{\pm}, \i  \tilde{H}])P^{+}.
\end{array}
 \]
 Since the second term in the rhs above is the adjoint of the first,  and since \eqref{e5.1} for $P^{+}$ implies \eqref{e5.1} for $P^{-}$, it suffices to solve
 \begin{equation}
 \label{e5.5}
 P^{+}(\p_{t}P^{+}+ [P^{+}, \i \tilde{H}])P^{-}=0,
 \end{equation}
modulo  a smoothing error. Using \eqref{e5.3b} we obtain since $[P^{\pm}, H]=0$:
 \[
 \begin{array}{rl}
 &P^{+}\left(\p_{t}P^{+}+ [P^{+}, \i \tilde{H}]\right)P^{-}\\[2mm]
 =&P^{+}\left(\p_{t}P^{+}+ P^{+}HP^{+}S- S P^{-}H P^{-}+ F_{-1}(S)\right)P^{-}
 \end{array}
 \]
 We use now \eqref{e5.-0} denoting  the scalar operator $\epsilon(t, \rx, D_{x})$ simply by $\epsilon$ and obtain:
 \[
 \begin{array}{rl}
& P^{+}HP^{+}S- S P^{-}H P^{-}= \epsilon S+ S \epsilon+ R_{0}^{+}S- S R_{0}^{-}\\[2mm]
=&2 \epsilon S+ [S, \epsilon]+ R_{0}^{+}S- S R_{0}^{-}.
 \end{array}
 \]
 The maps $S\mapsto R_{0}^{+}S$, $S\mapsto SR_{0}^{-}$ belong to $\cF_{-1}$ by Lemma \ref{lemma5.2b}, and so does  the map $S\mapsto [\epsilon, S]$, since $\epsilon$ is scalar. 

Therefore the equation \eqref{e5.5} can be rewritten as
\[
\p_{t}P^{+}+ 2\epsilon S+ F_{-1}(S)=0,
\]
or equivalently as
\begin{equation}
\label{e5.6}
S= - (2\epsilon)^{-1}\p_{t}P^{+}+ F_{-2}(S),
\end{equation}
where $F_{-2}: S\mapsto - (2 \epsilon)^{-1}F_{-1}(S)$ belongs to $\cF_{-2}$. We apply Prop. \ref{prop5.1} to solve \eqref{e5.6}. We find $S\in \cinfb(I; \Psi^{-1})$, unique modulo  $\cinfb(I; \cW^{-\infty})$ such that
\[
\p_{t}P^{+}+ 2 \epsilon S+ F_{-1}(S)\in \cinfb(I; \cW^{-\infty}),
\]
and hence
\[
\p_{t}P^{+}+ [P^{+}, \i \tilde{H}]= R_{-\infty}\in \cinfb(I; \cW^{-\infty}).
\]
Differentiating $ \tilde{\cU}(s,t)P^{+}(t) \tilde{\cU}(t,s)$ w.r.t. $t$ and using again Lemma \ref{lemma5.2} we obtain that $ \tilde{\cU}(t,s)P^{\pm}(s)- P^{\pm}(t)\tilde{\cU}(t,s)\in \cinfb(I^{2}; \cW^{-\infty})$ and hence if $ \tilde{P}^{\pm}(t)= \e^{- \i R(t)}P^{\pm}(t)\e^{\i R(t)}$ we have:
\[
\cU(t, s)\tilde{P}^{\pm}(s)-\tilde{P}^{\pm}(t)\cU(t,s)\in \cinfb(I^{2}; \cW^{-\infty}).
\]
Clearly $\tilde{P}^{+}(t)+ \tilde{P}^{-}(t)=\one$ and since $R(t)=R^{*}(t)$, $\e^{\i R(t)}$ is unitary and hence $\tilde{P}^{\pm}$ are selfadjoint, so we have proved (1) and (2) of the proposition. 

It remains to check (3). 
We set  $Q_{\pm}(t, \rx, \p_{t}, \p_{\rx})= \p_{t}\mp \i \epsilon(t)+ R^{\pm}_{0}(t)$, considered as an operator acting on $M\times \Sigma$ and denote by $A(t, \rx, \rx')\in \cD'(M\times \Sigma; L(S(\Sigma)))$ the distributional kernel of $\cU(\cdot, s)\tilde{P}^{\pm}(s)$.  Then $Q_{\pm}A\in \cinf(M\times \Sigma; L(S(\Sigma)))$. Since $Q_{\pm}$ is not a classical pseudodifferential operator on $M\times \Sigma$, we cannot directly apply the microlocal regularity to obtain (3).
Instead we use an argument from \cite[Lemma 6.5.5]{DH}. We fix a scalar pseudodifferential operator $Q_{0}\in \Psi^{0}_{\rm ph}(M\times \Sigma)$ with principal symbol $\chi(\frac{|\tau|+ |k'|}{|k|})$, where $\chi\in \coinf(\rr)$ is equal to $1$ on $[-C, C]$. Then $Q_{0}Q_{\pm}$ is a classical pseudodifferential operator  of order $1$ on $M\times \Sigma$, with principal symbol \[
 \i\chi(\frac{|\tau|+ |k'|}{|k|})(\tau\mp( k\dual h^{-1}(t, \rx)k)^{\12}).
 \]
   Since $Q_{0}Q_{\pm}A$ is smooth and $Q_{0}Q_{\pm}$ is elliptic in $\{|\tau|+ |k'|\leq C|k|, \ \tau\mp (k\cdot h^{-1}(t, \rx)k)^{\12}\neq 0\}$ we obtain taking $C$ arbitrarily large that
   \[
   \WF(A)'\subset (\cN^{\pm}\cup \cF)\times T^{*}\Sigma,
   \]
   as claimed. \qed

\subsection{Pure Hadamard states}
We consider now  a spacetime $(M, g)$ and a Dirac operator $D$  satisfying the hypotheses ${\rm (H)}$, ${\rm (M)}$ in Subsect. \ref{sec5.0}.

We will set  
\[
c^{\pm}\defeq \e^{\frac{1-n}{2} u}\tilde{P}^{\pm}(0)\e^{\frac{n-1}{2}u},
\] where the projections $\tilde{P}^{\pm}(t)$ are constructed in Prop. \ref{prop5.2}. We recall that we identified $\Sigma_{0}= \{0\}\times \Sigma$ with $\Sigma$ and that we use the Hermitian forms
\[
(\psi| \psi)_{M}= \int_{M}\bar{\psi}\dual\beta \psi dVol_{g},
\ (f| f)_{\Sigma}= \int_{\Sigma}\bar{f}\dual\beta f dVol_{h}
\]
to identify sesquilinear forms with linear operators.  We recall also that $\rho_{t}$ is the trace on $\Sigma_{t}$, $U_{\Sigma}f$ is the unique solution of the Cauchy problem:
\[
\left\{  
\begin{array}{l}
D\psi= 0,\\
\rho_{0}\psi= f
\end{array},\right. \ f\in \coinf(\Sigma, S(\Sigma)),
\]
and  $U(t,s)$ the Cauchy evolution operator for $D$, i.e. 
\beq\label{defdecauchy}
U(t,s)f\defeq \rho_{t}\psi, 
\hbox{ where }
\left\{  
\begin{array}{l}
D\psi= 0,\\
\rho_{s}\psi= f.
\end{array}, \ f\in \coinf(\Sigma_{s}, S(\Sigma)).
\right.
\eeq
\begin{theoreme}\label{maintheorem}
The maps
\[
\lambda_{\Sigma}^{\pm}= \i \gamma(n) c^{\pm}
\]
are the Cauchy surface covariances of a pure Hadamard state $\omega$ for $D$ on $M$.
\end{theoreme}
\proof 
We combine the results recalled in Subsect. \ref{sec1b.4} on conformal transformations  with Prop. \ref{prop5.2} and Prop. \ref{prop15.0a}. We obtain that $\lambda^{\pm}$ are the Cauchy surface covariances of a pure Hadamard state for $D$. \qed

\subsection{Spacetime covariances and Feynman inverses}
We now describe the 'time kernels' of the spacetime covariances $\Lambda^{\pm}$ and Feynman inverse $G_{\rm F}$ associated to the Hadamard state $\omega$ constructed in Thm. \ref{maintheorem}. For simplicity we assume that  we are in the model case considered in Subsect. \ref{sec5.1}.
Formulas in the general case can easily be obtained from  the identities in Subsect. \ref{sec1b.4}.

We recall that in the model case we have  $c^{\pm}= \tilde{P}^{\pm}(0)$.
\begin{proposition}\label{lastprop}
  The spacetime covariances $\Lambda^{\pm}$ and Feynman inverse $G_{\rm F}$ of the state $\omega$ in Thm. \ref{maintheorem} are given by:
\[
\begin{array}{l}
\Lambda^{\pm}v(t)= \int_{I} \Lambda^{\pm}(t, s)v(s)ds,\\[2mm]
 G_{\rm F}v(t)= \int_{I}G_{\rm F}(t, s) v(s)ds, \ v\in \coinf(I; \coinf(\Sigma;S(\Sigma))),
\end{array}
\]
for
\beq\label{e5.9}
\begin{array}{rl}
i)&\Lambda^{\pm}(t,s)= \i U(t, 0)c^{\pm}U(0,s )\gamma(e_{0}),\\[2mm]
ii)&G_{\rm F}(t, s)= U(t, 0)\left( \theta(t-s)c^{+}- \theta(s-t)c^{-}\right)U(0, s)\gamma(e_{0}).
\end{array}
\eeq

\end{proposition}
\proof We have:
\[
\Lambda^{\pm}= (\rho_{0}G)^{*}\i \gamma(e_{0})c^{\pm}\rho_{0}G = \i U_{\Sigma}c^{\pm}\rho_{0}G.
\]
The retarded/advanced inverses $\cG_{\rm ret/adv}$ for $\cD= \p_{t}- \i H(t)$ are given by
\[
\cG_{\rm ret/adv}v(t)= \int_{I}\cG_{\rm ret/adv}(t, s) v(s)ds,
\]
for 
\[
\cG_{\rm ret/adv}(t,s): \pm\theta(\pm(t-s))\cU(t,s),
\]
$\theta= \one_{\rr^{+}}$ being the Heaviside function. 
From \eqref{e3.8b} we obtain that
\begin{equation}
\label{e5.7}
U(t,s)= |h_{t}|^{-\frac{1}{4}}|h_{0}|^{\frac{1}{4}}
\cT(t, 0)\cU(t,s)\cT(0,s)|h_{s}|^{\frac{1}{4}}|h_{0}|^{-\frac{1}{4}}.
\end{equation}

 It follows that the retarded/advanced inverses for $D$ are given by
\[
G_{\rm ret/adv}v(t)= \int_{I}G_{\rm ret/adv}(t, s) v(s)ds,
\]
for 
\beq\label{e5.10}
G_{\rm ret/adv}(t,s): \pm\theta(\pm(t-s))U(t,s)\gamma(e_{0}),
\eeq
hence  the 'time kernel' $G(t,s)$ of $G= G_{\rm ret}- G_{\rm adv}$ equals
\[
G(t,s)= U(t,s)\gamma(e_{0}),
\]
which implies  the first statement in \eqref{e5.9}. The second follows then from 
\[
G_{\rm F}= \i^{-1} \Lambda^{+}+ G_{\rm adv}=  -\i^{-1} \Lambda^{-}+ G_{\rm ret},
\]
 and \eqref{e5.10}. \qed
 
\section{Hadamard states on arbitrary spacetimes}\label{sec7}\init
In this section we give two proofs of existence of Hadamard states for Dirac fields on general globally hyperbolic spacetimes with a spin structure. The first is by the usual deformation argument and is related to the one given in \cite{MV}. The second uses a partition of unity and is due to \cite{GW1} for Klein-Gordon fields.
\subsection{The vacuum state associated to a Killing field}
We first recall the definition of the vacuum state for Dirac fields associated to a Killing vector field, due to \cite{DHo}.
On a Lorentzian manifold with a spin structure, the Lie derivative of a spinor field is defined as (see \cite{Kos}) :
\beq\label{kosmano}
\begin{array}{l}
\cL_{X}\psi= \nabla_{X}^{S}\psi+ \frac{1}{8}((\nabla_{a}X)_{b}- (\nabla_{b}X)_{a})\gamma^{a}\gamma^{b}\psi,\\[2mm]
 \psi\in \cinf(M; S(M)), \ X\in \cinf(M; TM).
\end{array}
\eeq
The Lie derivative extends in the obvious way to  $\cinf(M; L(S(M)))$ by setting
\[
\cL_{X}A\psi\eqdef (\cL_{X}A)\psi+ A \cL_{X}\psi.
\]
If  $X$ is a complete Killing vector field and the mass term $m$ in Def. \ref{def1.2} satisfies $\cL_{X}m=0$, then $[D, \cL_{X}]=0$. It follows that the flow $\phi_{s}$ generated by $\cL_{X}$ preserves $\Sol(D)$. One can easily show using  \ref{e10.4} that $\phi_{s}$  preserves the Hilbertian scalar product $\nu$ defined in \ref{defidef}.

It hence defines a unique strongly continuous unitary group $(\e^{\i sH})_{s\in \rr}$ on the completion of $(\Sol(D), \nu)$ whose generator $H$ is by Nelson's invariant domain theorem equal to the closure of $\i^{-1}\cL_{X}$ on $\Sol(D)$.

Alternatively we can unitarily identify (the completions of) $(\Sol(D), \nu)$ and $(\coinf(\Sigma; S(\Sigma)), \nu_{\Sigma})$ by $\rho_{\Sigma}$ as in \eqref{unitarymaps}. The image of $H$ by this identification is denoted by $H_{\Sigma}$ acting on $L^{2}(\Sigma; S(\Sigma))$. It is essentially selfadjoint on $\coinf(\Sigma; S(\Sigma))$.
\begin{definition}\label{defdevide}
 Let $(M, g)$ be  globally hyperbolic spacetime with a complete Killing field $X$ and $D= \slashed{D}+ m$ a Dirac operator on $(M, g)$ with $\cL_{X}m=0$. 
 Assume moreover that
 \begin{equation}
 \label{kerokero}
 \Ker H_{\Sigma}= \{0\}.
 \end{equation}
 Then  the {\em vacuum state} $\omega^{\rm vac}$ associated to $X$ is the quasi-free state defined by the Cauchy surface covariances:
 \[
 \lambda^{\pm{\rm vac}}\defeq \i \gamma(n)\one_{\rr^{\pm}}(H_{\Sigma}).
 \]
 \end{definition}
Unlike  the bosonic case, $X$ does not need to be time-like in order to be able to define the associated vacuum state.

Let us now compute $H_{\Sigma}$ in the ultrastatic case, i.e. if $M= \rr\times \Sigma$, $g= -dt^{2}+ h(\rx)d\rx^{2}$ and the mass term $m$ in Def. \ref{def1.2} is a constant scalar.  The restriction of the spinor bundle $S(M)$ to $\{t\}\times \Sigma$ is then independent on $t$ and denoted by $S(\Sigma)$.
The flow $\phi_{s}$ associated to $\cL_{\p_{t}}$ is simply the time translations $\phi_{s}\psi(t, \cdot)= \psi(t-s, \cdot)$, and  setting $\gamma_{0}= \gamma(\p_{t})$ we have
 \[
 D= - \gamma_{0}(\p_{t} -\i H_{\Sigma}), 
 \]
 and \[
 H_{\Sigma}= \i \gamma_{0}(\gamma(e_{a})\nabla^{S}_{e_{a}}+ m)\eqdef H_{0\Sigma}+ \i \gamma_{0}m,
 \]
 where $(e_{a})_{1\leq a\leq d}$ is a local orthonormal frame of $T\Sigma$ for $h$.  Since $H_{0\Sigma}\gamma_{0}= - \gamma_{0}H_{0\Sigma}$ we obtain $H_{\Sigma}^{2}= H_{0\Sigma}^{2}+ m^{2}$, hence  $\Ker H_{\Sigma}= \{0\}$ if $m>0$.
 
 By \cite[Thm. 5.1]{SV2} we know that  the vacuum state  $\omega^{\rm vac}$ is a Hadamard state.  Let us sketch a direct proof using pseudodifferential calculus.   Let $\Psi^{m}_{\rm c}(\Sigma; S(\Sigma))$ the space  of properly supported supported classical pseudodifferential operators on $\Sigma$ and $\Psi^{m}(\Sigma; S(\Sigma))=  \Psi^{m}_{\rm c}(\Sigma; S(\Sigma))+ W_{-\infty}$, where $W_{-\infty}$ is the space of smoothing operators  acting on $S(\Sigma)$. 
 
 Since $H_{\Sigma}$ is elliptic in $\Psi^{1}(\Sigma; S(\Sigma))$, selfadjoint with $0\not\in \sigma(H_{\Sigma})$, we obtain by the same arguments as in the proof of Prop. \ref{prop4.3} that $\one_{\rr^{\pm}}(H_{\Sigma})\in \Psi^{0}(\Sigma; S(\Sigma))$. From Clifford relations we obtain that $\sigma_{\rm pr}(H_{\Sigma})^{2}(\rx, k)= k\dual h^{-1}(\rx)k\one$, which implies that 
 \[
 H_{\Sigma}\one_{\rr^{\pm}}(H_{\Sigma})= (\pm \epsilon(\rx, D_{\rx})+ R_{0}^{\pm}(\rx, D_{\rx})\one_{\rr^{\pm}}(H_{\Sigma}),
 \]
  where $R_{0}^{\pm}\in \Psi^{0}$ and $\epsilon$ is a scalar pseudodifferential operator in $\Psi^{1}(\Sigma; S(\Sigma))$ with principal symbol $(k\dual h^{-1}(\rx)k)^{\12}$.
  
 Since $U_{\Sigma}f(t, \cdot)=\e^{\i tH_{\Sigma}}f(\cdot)$,  we obtain that $(\p_{t}\mp \i \epsilon+ R_{0}^{\pm})U_{\Sigma}c^{\pm{\rm vac}}=0$.
   We argue then as in the proof of Prop. \ref{prop5.2} (3) to obtain that
 \[
 \WF(U_{\Sigma}c^{\pm{\rm vac}})'\subset(\cN^{\pm}\cup \cF)\times T^{*}\Sigma, \ \cF= \{\xi=0\}\subset  T^{*}M.
 \]
 By Prop. \ref{prop15.0a} this implies that $\omega_{\rm vac}$ is a Hadamard state.

\subsection{Deformation of Dirac operators}\label{sec7.2}
Let $(M, g)$ be a globally hyperbolic spacetime and $\Sigma$ a smooth space-like Cauchy surface. By the Bernal-Sanchez theorem \cite{BS1, BS2} we can assume that 
\[
M= \rr_{t}\times \Sigma_{\rx}, \ g= - c^{2}(t, \rx)dt^{2}+ h(t, \rx)d\rx^{2},
\]
where $c\in \cinf(M, \rr)$, $c>0$, $\rr\ni t\mapsto h(t, \rx)d\rx^{2}$ is a family of Riemannian metrics on $\Sigma$ and 
\[
\Sigma_{s}=\{s\}\times \Sigma
\]
 are Cauchy surfaces for $g$.

 After a conformal transformation we can assume that $c\equiv 1$. 
 We  fix some complete Riemannian metric $h_{\rm us}(\rx)d\rx^{2}$ on $\Sigma$  and set
 \[
 g_{\rm us}= -dt^{2}+ h_{\rm us}(\rx)d\rx^{2}.
 \]
 The existence of a  {\em globally hyperbolic }interpolating metric $g_{\rm int}$ such that 
 \[
 g_{\rm int}= \left\{\begin{array}{l}
 g\hbox{ in }t> 1,\\
 g_{\rm us}\hbox{ in }t<1 
 \end{array}\right.
 \]
 was apparently taken for granted. It was realized only recently by Sanchez in \cite{Sa} that the simple constructions of $g_{\rm int}$ appearing in the literature may fail to be globally hyperbolic. However Müller has proved in \cite[Thm. 3]{Mu} that globally hyperbolic  interpolating metrics $g_{\rm int}$ exist. 
 
  Let us now assume that $(M, g)$ has a spin structure. Since $M$ is orientable,  this spin structure is unique.   
 
 We saw in Subsect. \ref{spinprod} that  it  induces a  spin structure $P\Spin(\Sigma, h_{0})$ on $\Sigma$.  From this spin structure on $(\Sigma, h_{0})$ we obtain as in Subsect. \ref{spinprod} a spin structure on $(M, g_{\rm int})$.  We also saw that if $S(M)$ is the associated spinor bundle, its restriction to $\{t\}\times \Sigma$ is independent on $t$ and denoted by $S(\Sigma)$. Over $\{t>1\}$ this spin structure   coincide with the original spin structure on $(M, g)$. The same is true for the associated spinor bundles.
 
 Let now $D= \slashed{D}+ m$ a Dirac operator on $(M, g)$ as in Subsect. \ref{sec1.4}. We fix a strictly positive constant $m_{\rm us}>0$ and an interpolating mass $m_{\rm int}$ such that
 \[
 m_{\rm int}(t, \rx)= \left\{\begin{array}{l}
 m(t, \rx )\hbox{ for }t>1, \\[2mm]
 m_{\rm us}\one\hbox{ for }t< -1.
 \end{array}  \right. 
\]
We denote by $D_{\rm int}= \slashed{D}+ m_{\rm int}$, resp.  $D_{\rm us}= \slashed{D}+ m_{\rm us}$  the Dirac operator obtained from the spin structure on $(M, g_{\rm int})$ resp. on $(M, g_{\rm us})$.  We have
\[
 D_{\rm int}=\left\{\begin{array}{l}
 D\hbox{ for }t>1, \\[2mm]
 D_{\rm us}\hbox{ for }t< -1.
 \end{array}  \right. 
\]
We denote by $U_{\rm int}(t, s)$ the Cauchy evolution operator for $D_{\rm int}$ see \eqref{defdecauchy}.
\subsection{Existence of Hadamard states}\label{sec7.3}
We now sketch a proof of existence of Hadamard states for Dirac fields on an arbitrary globally hyperbolic spacetime with a spin structure.
\begin{theoreme}\label{thm7.1}
 Let $(M, g)$ an even dimensional globally hyperbolic spacetime with a spin structure and $D$ a Dirac operator on $(M, g)$. Then there exist pure Hadamard states for $D$.
\end{theoreme}
\proof 
The proof is completely analogous to the case of scalar bosonic fields.
We reduce ourselves to the situation at the beginning of Subsect. \ref{sec7.2} by the Bernal-Sanchez theorem and  a conformal transformation. 
Let $\lambda^{\pm{\rm vac}}$ be the Cauchy surface covariances of the vacuum state $\omega^{\rm vac}$ for $D_{\rm us}$, on the Cauchy surface $\Sigma_{-2}$.
Since $D_{\rm int}= D_{\rm us}$ in $\{t<-1\}$,  $\lambda^{\pm{\rm vac}}$  satisfy \eqref{troup.100} in Prop. \ref{prop15.0a} for $D_{\rm int}$ over $\{t<-1\}$. By propagation of singularities, (see the end of the proof of Prop. \ref{prop15.0a} for the precise argument), we obtain that 
$U_{\rm int}(-2, 2)^{*}\lambda^{\pm{\rm vac}}U_{\rm int}(-2, 2)$ satisfy \eqref{troup.100} near $\Sigma_{2}$ for $D_{\rm int}$ and hence for $D$.
Therefore $U_{\rm int}(-2, 2)^{*}\lambda^{\pm{\rm vac}}U_{\rm int}(-2, 2)$ are the Cauchy surface covariances on $\Sigma_{2}$ of a Hadamard state $\omega$ for $D$. Since $\omega^{\rm vac}$ is a pure state, we see that $\omega$ is pure, by Prop. \ref{hurlub}.  \qed

\subsubsection{An alternative construction}
Let us finally give an alternative construction of Hadamard states on general spacetimes, analogous to a construction in \cite{GW1} for Klein-Gordon fields.  Unlike the previous construction it does not produce pure states. 

\def\tV{\tilde{V}}
 We fix a smooth spacelike Cauchy surface $\Sigma$ in $(M, g)$ and we can assume that $M= \rr\times \Sigma$.  We identify $\Sigma$ with $\{0\}\times \Sigma$. We fix a causally compatible neighborhood $U$ of $\Sigma$ in $M$, an atlas $(\tV_{i}, \chi_{i})_{i\in \nn}$ of $\Sigma$  with $\chi_{i}: \tV_{i}\tosim B_{d}(0, 1)$  and $\tV_{i}$ precompact, open  sets $V_{i}\subset \tV_{i}$ and constants $\delta_{i}>0$ such that
 \begin{equation}
 \label{blurbi}
 \begin{array}{rl}
 i)& \cup_{i\in \nn} V_{i}= \Sigma,\\[2mm]
 ii)& x\in U, J(x)\cap V_{i}\neq \emptyset\Rightarrow y\in ]-\delta_{i}, \delta_{i}[\times \tV_{i}\eqdef U_{i}.
 \end{array}
 \end{equation}
  We transport by $\phi_{i}: (t, \rx)\mapsto (t,\chi_{i}(\rx))$  the metric, spin structure and the Dirac operator $D$ over $U_{i}$ to $\phi_{i}(U_{i})$. We can extend the spin structure over $\phi_{i}(U_{i})$ to $\rr^{1+d}$ and obtain a Dirac operator $D_{i}$ on $\rr^{1+d}$ for a metric and spin structure which are of bounded geometry for the 'uniform' reference metric $dt^{2}+ dx^{2}$. 
  
  By Sect. \ref{sec5} we can construct a Hadamard state $\omega_{i}$ for $D_{i}$. Let $\lambda^{\pm}_{i}$ its Cauchy surface covariances on $\Sigma$ and $1= \sum_{i}u_{i}^{2}$ a partition of unity on $\Sigma$ subordinate to the covering $(V_{i})_{i\in \nn}$. We set
  \[
  \lambda^{\pm}= \sum_{i\in \nn}u_{i}^{*}\circ (\varphi_{i}^{-1})^{*}\lambda_{i}^{\pm}\circ u_{i}.
  \]
  Since $\nu_{\Sigma}= \sum_{i\in \nn}u_{i}^{*}\circ \nu_{\Sigma}\circ u_{i}$ the conditions \eqref{plic} and \eqref{troup.100} are clearly satisfied by $\lambda^{\pm}$. Therefore we obtain a Hadamard state for $D$ on $(M, g)$.

\appendix
  \section{}\label{app}\init

 \subsection{}\label{app.1}
 Let $\pmb{\eta}= {\rm diag}(-1, 1, \dots, 1)$ the Minkowski metric on $\rr^{1, d}$, $n= 1+d$. We set
 \[
 \begin{array}{l}
T_{n}^{<}\defeq \{\pmb{t}\in M_{n}(\rr): \pmb{t}_{ij}=0, \ j>i\}, \\[2mm]
T_{n}^{>} \defeq \{\pmb{t}\in M_{n}(\rr): \pmb{t}_{ij}=0, \ j<i\}, \\[2mm]
S_{n}\defeq  \{\pmb{s}\in M_{n}(\rr): \pmb{s}_{ij}= \pmb{s}_{ji}\},
\end{array}
\]
 which are $\rr$-vector spaces of dimension $n(n+1)/2$.
 
\begin{lemma}\label{lemma-app.1}
 There exists $V_{0}$ neighborhood of $\eta$ in $S_{n}$ and $U_{0}$ neighborhood of $\one_{n}$ in $T_{n}^{>}$ and $F: V_{0}\to U_{0}$ a smooth diffeomorphism such that
 \[
F(\pmb{g}) \pmb{g} ^{t}\!F(\pmb{g})= \pmb{\eta}, \ \forall \pmb{g}\in V_{0}.
\]
\end{lemma}
 \proof 
 Let $H: T_{n}^{<}\ni \pmb{t}\mapsto\pmb{t}\pmb{\eta}^{t}\! \pmb{t}\in S_{n}$. The differential of $H$ at $\one_{n}$ is
 \[
DH_{\one}: \pmb{t}\mapsto \pmb{t}\pmb{\eta}+ \pmb{\eta}^{t}\! \pmb{t}.
\]
 The equation  $ \pmb{t}\pmb{\eta}+ \pmb{\eta}^{t}\! \pmb{t}=\pmb{s}$ for $\pmb{s}\in S_{n}$ is written as
 \[
\pmb{t}_{ij}\pmb{\eta}_{jj}+ \pmb{\eta}_{ii}\pmb{t}_{ji}= \pmb{s}_{ij}, \  \pmb{\eta}_{11}= -1, \ \pmb{\eta}_{ii}= 1\ j>1,
\]
which is solved by
\[
\pmb{t}_{ij}= \pmb{\eta}_{jj}^{-1}\pmb{s}_{ij}, \hbox{ for }j<i, \ \pmb{t}_{ii}= ( 2\pmb{\eta}_{ii})^{-1}\pmb{s}_{ii}, \ \pmb{t}_{ij}=0\hbox{ for }j>i.
\]
It follows that $DH_{\one}: T_{n}^{<}\to S_{n}$ is surjective hence bijective. By the local inversion theorem, there exists $W_{0}$ neighborhood of $\one_{n}$ in $T_{n}^{<}$, $V_{0}$ neighborhood of $\pmb{\eta}$  in $S_{n}$ such that $H: W_{0}\tosim V_{0}$ is a smooth diffeomorphism. We compose $H$ with the smooth diffeomorphism $W_{0}\ni \pmb{g}\mapsto \pmb{g}^{-1}\in U_{0}$, where $U_{0}$ is a neighborhood of $\one_{n}$ in $T_{n}^{>}$ to obtain the lemma. \qed
\subsection{}\label{app.2}
  \begin{proposition}\label{propadd1.1}
Let $G\to P\xrightarrow{\pi}M$  and $G'\to P'\xrightarrow{\pi'}M$ two principal bundles over $M$ with structure groups $G, G'$. Let $\varphi: G'\to G$ a group morphism and  a local proper diffeomorphism and   $\chi: P'\to P$  a  bundle morphism and a local diffeomorphism. Assume that the following diagram commute:
\beq\label{e1.-1}
 \begin{tikzcd}
G' \arrow[r] \arrow[dd, "\varphi"] & P' \arrow[rd, "\pi'"] \arrow[dd, "\chi"] &   \\
                   &                         & M \\
G \arrow[r]            & P \arrow[ru, "\pi"]            &  
\end{tikzcd}
\eeq
ie $\chi$ is a morphism of principal bundles over $M$.

 Then  modulo an isomorphism of $G'$-principal bundles over $M$ one can assume  $G$ and $G'$ are locally trivialized over a common covering $(U_{i})_{i\in I}$ of $M$ and that
 \[
 t_{ij}(x)=  \varphi(t'_{ij}(x)), \ x\in  U_{i}\cap U_{j},
 \]
 where $t_{ij}: U_{i}\cap U_{j}\to G$, $t'_{ij}: U_{i}\cap U_{j}\to G'$ are the transition maps of $P$, $P'$.
  \end{proposition}
\proof  
Let  $(U_{i}^{(\prime)})_{i\in I^{(\prime)}}$   a bundle atlas for $P^{(\prime)}$,  $\Psi_{i}^{(\prime)}: (\pi^{(\prime)})^{-1}(U_{i}^{(\prime)})\to U_{i}^{(\prime)}\to G^{(\prime)}$ the local trivializations of $P^{(\prime)}$, $\Phi_{i}^{(\prime)}= (\Psi_{i}^{(\prime)})^{-1}$.  Without loss of generality we can assume that $U_{i}^{(\prime)}$ are simply connected.

If we set $m_{i}^{(\prime)}(x)= \Phi_{i}^{(\prime)}(x, e^{(\prime)})$, where $e^{\prime}$ is the identity in $G^{(\prime)}$, then  since $G^{(\prime)}$ acts transitively on the fibers of $P^{(\prime)}$, we have
\beq\label{e1.-1a}
\Phi_{i}^{(\prime)}(x, g^{(\prime)})= m_{i}^{(\prime)}(x)\cdot g^{(\prime)}, \ x\in U_{i}^{(\prime)}, \ g^{(\prime)}\in G^{(\prime)}.
\eeq

Since $U_{i}$ is simply connected,  there exists $\tilde{m}_{i}': U_{i}\to \pi^{\prime -1}(U_{i})$ such that $m_{i}= \chi\circ \tilde{m}_{i}'$, i.e. a lift of $m_{i}$ through $\chi$.  Since $G'$ acts transitively on the fibers of $P'$, we can define:
\beq\label{e1.-2a}
\tilde{\Phi}'_{i}: U_{i}\times G'\ni 
(x, g')\mapsto \tilde{m}_{i}'(x)\cdot g'\in \pi'^{-1}(U_{i})
\eeq
Let us study some properties of $\tilde{\Phi}^{\prime}_{i}$. First 
\[
\begin{array}{rl}
&\pi'\circ \tilde{\Phi}_{i}'(x, g')= \pi'(\tilde{\Phi}_{i}^{\prime}(x)\cdot g')= \pi'\circ \tilde{m}_{i}'(x)\\[2mm]
=&\pi\circ \chi\circ \tilde{m}_{i}'(x)= \pi\circ m_{i}(x)= x.
\end{array}
\]
Next
\beq\label{e1.1}
\begin{array}{rl}
&\chi(\tilde{\Phi}_{i}'(x, g'))= \chi(\tilde{m}_{i}'(x)\cdot g')= \chi(\tilde{m}_{i}'(x))\cdot \varphi(g')= m_{i}(x)\cdot \varphi(g')\\[2mm]
=& \Phi_{i}(x, e)\cdot \varphi(g')= \Phi_{i}(x, \varphi(g')).
\end{array}
\eeq
Let us set $\tilde{\Psi}_{i}'= \tilde{\Phi}_{i}^{\prime -1}: \pi^{\prime -1}(U_{i})\tosim U_{i}\times G'$. The family of trivialisations $(\tilde{\Psi}_{i}')_{i\in I}$, together with the  right action of $G'$ on $P'$ and the base projection $\pi': P'\to M$ define hence a $G'$-principal bundle $G\to\tilde{P}'\xrightarrow{\pi}M$, (equal to $P'$ as a set).  If $\tilde{t}'_{ij}: U_{ij}\to G'$ are the associated transition maps we have
\beq\label{e1.1b}
\begin{array}{rl}
&(x, t_{ij}(x)\cdot \varphi(g'))= \Psi_{i}\circ \Phi_{j}(x, \varphi(g'))= \Psi_{i}\circ \chi\circ \tilde{\Phi}_{j}^{\prime}(x, g')\\[2mm]
=& \Psi_{i}\circ \chi\circ \tilde{\Phi}_{i}^{\prime}(x, t'_{ij}(x)g')= (x, \varphi(\tilde{t}_{ij}')(x)\varphi(g')),
\end{array}
\eeq
hence $t_{ij}= \varphi(\tilde{t}'_{ij})$. Let now $G'\to \tilde{P}'\xrightarrow{\pi'}M$ the $G'$-principal bundle associated with base projection $\pi'$ and trivializations $\tilde{\Psi}_{i}'= (\tilde{\Phi}_{i}')^{-1}: (\pi')^{-1}(U_{i})\to U_{i}\times G'$.

 It remains to prove that $\tilde{P}'$ and $P'$ are isomorphic as $G'$-principal bundles.   For $\pmb{i}= (i, i')\in \pmb{I}= I\times I'$ we set $U_{\pmb{i}}= U_{i}\cap U'_{i'}$, so that $(U_{\pmb{i}})_{\pmb{i}\in \pmb{I}}$ is a common bundle atlas for $P'$ and $\tilde{P}'$.

 Since $G'$ acts transitively on the fibers of $P'$ and $\tilde{P}'$ we can find maps $\chi_{\pmb{i}}: U_{\pmb{i}}\to G'$ such that
 \[
 m_{i}'(x)= \tilde{m}'_{i}(x)\cdot \chi_{\pmb{i}}(x), x\in U_{\pmb{i}}, \ \pmb{i}= (i, i'),
 \]
 which by \eqref{e1.-1a}, \eqref{e1.-2a}   implies that
 $\Phi_{i'}'(x, g')= \tilde{\Phi}_{i}'(x, g')\cdot \chi_{\pmb{i}}(x)$, $x\in  U_{\pmb{i}}$ and hence that $\tilde{P}'$ and $P'$ are isomorphic.  \qed

\end{document}